\documentclass[12pt]{article}
\usepackage{amsmath}
\usepackage{graphicx,psfrag,epsf}
\usepackage{enumerate}
\usepackage{natbib}
\usepackage{url} 
\usepackage{xcolor}

\newcommand{\blind}{0}

\RequirePackage{amsthm,amsmath,amsfonts,amssymb}
\theoremstyle{remark}
\theoremstyle{plain}

\newtheorem{proposition}{Proposition}[section]
\newtheorem{remark}{Remark}[section]

\begin{document}

\def\spacingset#1{\renewcommand{\baselinestretch}%
{#1}\small\normalsize} \spacingset{1}


\if0\blind
{
  \title{\bf  A Characterization of Most(More) Powerful Test Statistics with Simple Nonparametric Applications}
  \author{Albert Vexler 
    \hspace{.2cm}\\
    Department of Biostatistics, The State University of New York at Buffalo, 
    \\
    Buffalo, NY 
    \\
    and \\
    Alan D. Hutson \\
    Department of Biostatistics and Bioinformatics, \\
    Roswell Park Comprehensive Cancer Center, Buffalo, NY
}
  \maketitle
} \fi

\if1\blind
{
  \bigskip
  \bigskip
  \bigskip
  \begin{center}
    {\LARGE\bf Title}
\end{center}
  \medskip
} \fi

\bigskip
\begin{abstract}
Data-driven most powerful tests are statistical hypothesis decision-making 
tools
that deliver the greatest power against a fixed null hypothesis among all corresponding data-based tests of a given size. When the underlying data distributions are known, the likelihood ratio principle 
can be applied to conduct
 most powerful tests. 
  Reversing 
 this notion, we consider the following questions. (a) Assuming a test statistic, say $T$, is given, how can we transform $T$ to improve the power of the test?
(b) 
Can $T$ be used to generate the most powerful test?
 (c) How does one compare test statistics with respect to an attribute of the desired most powerful decision-making procedure? 
To examine these questions, we propose  one-to-one mapping of the term "most powerful" to the distribution properties of a given test statistic via 
 matching characterization.
This form of characterization has practical applicability and aligns well  with the general principle of sufficiency.
Findings indicate that to improve a given test, we can employ relevant ancillary statistics that do not have changes in their distributions with respect to tested hypotheses.
As an example, the present method is illustrated by modifying the usual t-test under nonparametric settings.
%
%
Numerical studies based on generated data and a real-data set confirm that the proposed approach can be useful in practice.
\end{abstract}

\noindent%
{\it Keywords:}  
Ancillary statistic; 
Likelihood ratio;
Most powerful test;
Nonparametric test;
Sufficiency;
t-test;
Test for median.
\vfill

\newpage
\spacingset{1.45} 
\section{Introduction}\label{sec1}
Methods for developing  and examining data-based  decision-making mechanisms  have been widely established in both  theoretical and experimental statistical frameworks.    
The common approach for evaluating  modern data-based testing algorithms follows the standards and foundations formulated nearly a century ago. In this context, for an extensive review and associated examples we refer the reader to~\cite{LehmanB}.
The criteria for which statistical tests are commonly competed against each other
uses the following prescription: 
1) type I error (TIE) rates of considered tests are fixed at the same level, say $\alpha$;
 and 2)  power levels of the tests are compared.   
This classical principle was largely created and advocated by J. Neyman and
E. S. Pearson in a series of substantive  papers published during 1928$-$1938 (e.g.,~\citealp{LehmanJASA}). In this framework, 
the likelihood methodology is associated with the likelihood
ratio concept, which allows for the development of powerful statistical inference
tools in decision-making tasks.

	In view of this, the likelihood ratio principle can be employed across a wide range of decision-making problems, although  likelihood ratio tests are not completely specified in many practical applications. Cases exist in which estimated parametric likelihood ratio statistics can have different formulations depending upon the underlying schemes of estimations.
	 Other times, the relevant likelihood functions may be quite complicated when, for example, the observations belong to correlated longitudinal data subject to some type of missing data mechanism.  There are other nonparametric scenarios that limit our ability to write corresponding likelihood ratio statistics. A systematic study of the inherent properties of likelihood ratios is necessary for proposing policies to advise on procedures for constructing, modifying, and/or selecting test statistics.

Without loss of generality, and in order to simplify the explanations of the main aim of this paper, we state the following formal notations. 
Assume we observe the underlying data $D$ with the goal of testing  a simple null hypothesis $H_0$ against its simple alternative $H_1$. To this end, we let $T=T(D)$  denote a real valued one-dimensional statistic based on $D$ that supports rejection of $H_0$, if $T(D)>C$, where $C$ is a fixed test-threshold. We write $\Pr_{k}$ and $\text{E}_{k}$ to denote the probability and expectation under $H_k$, $k\in\{0,1\}$, respectively.  
Throughout most of this paper, we suppose that the probability distributions $\Pr_{0}$, $\Pr_{1}$ of $D$ are absolutely continuous with respect to a given sigma finite measure $\tau$ defined over  $\Upsilon$, where $\Upsilon$ is an additive class of sets in a space, say $\mathbb{X}$, over which $D$ is distributed.  
Then, we have non-negative probability density functions $f_{0}$ and $f_{1}$ with respect to $\tau$ that satisfy    
\begin{eqnarray*}
	{\Pr}_k\left(D\in\upsilon\right)&=&\int_{\upsilon}f_k (x) d\tau(x),\quad \upsilon\in\Upsilon,\quad k\in\{0,1\}.
\end{eqnarray*}
Note that $f_0$ and $f_1$ need not belong to the same parametric family of
distributions.
Define the probability  density functions $f^{T}_0$ and $f^{T}_1$ of   $T$ such that
\begin{eqnarray*}
	{\Pr}_k\left\{T(D)\le t\right\}&=&\int_{\mathbb{X}}I\left(T(x)\le t\right)f_k (x) d\tau(x)\,=\,
	\int_{ {{\mathbb{R}}^{1}}} I\left(u\le t\right)f^{T}_k (u) du,
\end{eqnarray*}
  where $I(.)$ means the  indicator function, $H_k$ is assumed to be true, and $k\in\{0,1\}$. In this framework,  the likelihood ratio $\Lambda=\Lambda(D)=f_{1}(D)/f_{0}(D)$ is the most powerful (MP) test statistic, if $f_{0}(x)>0$ and $f_{1}(x)>0$, for all $x\in \mathbb{X}$. In scenarios when  the random observable data $D$ is multidimensional, $f_{0}$ and $f_1$ are generalized joint probability density with respect to $\tau$ (e.g.,~\citealp{LehmanAMS}).  

For the sake of simplicity, it will be assumed that in a case when a researcher plans to employ a test statistic $S(D)=L\left(\Lambda(D)\right)$, where $L(u)$ is a monotonically increasing function with the inverse function $W(L(u))=u$, we suppose the test statistic $T(D)=W(S(D))$ is in use and $T(D)$ is MP.   

The starting point of our study is associated with the following property of $\Lambda$ that can be found in~\cite{VexlerB} and is included  for 
the sake of completeness.
\begin{proposition}\label{pr1}
	The likelihood ratio statistic $\Lambda$ satisfies $f_{1 }^{\Lambda} (u)=u{\rm \; }f_{0}^{\Lambda} (u),$ 
	for all $u>0$. 
\end{proposition}
The proof is deferred to the online supplementary materials.

An interesting observation is that the likelihood ratio $f_{1 }^{\Lambda}/ f_{0}^{\Lambda}$, based on the likelihood ratio $\Lambda$  is itself, forms the likelihood ratio $\Lambda$, i.e., $f_{1 }^{\Lambda} (\Lambda)/f_{0}^{\Lambda} (\Lambda)=\Lambda$. Consider the situation when a value of a statistic or a single data point, say $X$, is observed.
The best transformation of $X$ for making a decision with respect to  $H_0$ against $H_1$ is the ratio $f_{1 }^{X} (X)/f_{0}^{X}(X)$. In the case $X=\Lambda$, the observed statistic cannot be improved, and the transformation $f_{1 }^{X} (X)/f_{0}^{X}(X)$ is invertible at $X=\Lambda$, since the likelihood ratio
$\Lambda$ is a root of the equation $f_{1 }^{X} (X)/f_{0}^{X}(X)=X$.


Let us for a moment assume  that we could improve Proposition~\ref{pr1} by including  the idiom "if and only if" in its statement, thereby asserting: a statistic $T$ has a likelihood ratio form if and only if $f_{1 }^{T} (u)=u{\rm \; }f_{0}^{T} (u)$. Then, since $\Lambda$ is the MP test statistic, having a given test statistic $T$, we will try to modify the structure of $T$, minimizing the distance between $f_{1 }^{T} (u)$ and $u{\rm \; }f_{0}^{T} (u)$, for at least some values of $u$. In this framework, comparing two test statistics, say $T$ and $B$, we can select $T$, if, for example, $\max_u\left|f_{1 }^{T} (u)-u{\rm \; }f_{0}^{T} (u)\right|$ $\le \max_u\left|f_{1 }^{B} (u)-u{\rm \; }f_{0}^{B} (u)\right|$.  
Note that in nonparametric settings we can approximate and/or estimate distribution functions of test statistics in many scenarios. Unfortunately, the simple statement of Proposition~\ref{pr1} cannot be used to characterize MP test statistics. For example, when $D=\left(X_1,X_2\right)$, where $X_1$ is from a normal distribution with $\text{E}_{0}\left(X_1\right)=0$, 
$\text{E}_{1}\left(X_1\right)=\mu$, and $\text{E}_{0}\left(X^2_1\right)=\text{E}_{1}\left(X_1-\mu\right)^2$, the statistic $T=f_{1 }\left(X_1\right)/f_{0}\left(X_1\right)$ satisfies 
$f_{1 }^{T} (u)=u{\rm \; }f_{0}^{T} (u),$ whereas the MP test statistic is $f_{1 }\left(X_1,X_2\right)/f_{0}\left(X_1,X_2\right)$. That is, to characterize MP tests, Proposition~\ref{pr1} needs to be modified.    
\begin{remark} {\em 
	In this article, we  avoid using the term "Uniformly Most Powerful" in order to be able to study cases when parameters do not play  an essential role in testing procedures as well as to consider situations where uniformly most powerful tests do not exist, e.g., nonparametric testing statements. 
	For example, let $H_0$ infer that observations are from a standard normal distribution vs. that the  observations follow a standard logistic distribution, under $H_1$. In this case, we have the likelihood ratio MP test,  whereas invoking the term "uniformly most powerful" can  be  misleading.  In Section~\ref{sc3}, the MP concept has a hypothetical context to which we aim to approach when we  develop  test procedures.  
	(See also item (iii) in Remark~\ref{rem2} and note (d) presented in Section~\ref{sc5}, in this aspect.) }
\end{remark}	
The goal of the present article is two-fold: (1) through an extension of Proposition~\ref{pr1}, we describe a way to characterize MP tests,
and (2) we exemplify the usefulness of the theoretical characterization of MP tests via  corrections of test statistics that can be  easily applied  in practice. 
Under this framework, Section~\ref{sc2} considers  one-to-one mapping of distribution properties of test statistics  to their ability to be most powerful. 
The proven  characterization is shown to be consistent with the principle of sufficiency in certain decision problems extensively  evaluated  by~\cite{Bah}.
 In Section~\ref{sc3}, to exemplify potential uses of the proposed  MP characterization, 
 we apply the theoretical concepts shown in Section~\ref{sc2} toward  demonstrating an efficient principle of improving commonly used test procedures via employing relevant ancillary statistics.
  Ancillary statistics have distributions that  do not depend on the competing hypotheses. 
   However,  we show  that ancillary statistics can make significant contributions to inference about the hypotheses of interest. 
For example, although it seems to be very difficult to compete against the well-known one sample  t-test for the mean, we assert that a simple modification of the t-test statistic can
increase its power.
This can be accomplished by accounting for the effect of population skewness on the distribution of the sample mean.
Section~\ref{sc3} demonstrates 
modifications of testing procedures that can be implemented under  nonparametric assumptions when there are no MP decision-making mechanisms. Then, in Section~\ref{sc4},  we show experimental evaluations that confirm high efficiency of the presented schemes in various situations.  Furthermore, as described in Section~\ref{data}, when used to analyze data from a biomarker study associated with myocardial infarction disease, the method proposed in Section~\ref{sc3} for one-sample testing about the median  is more sensitive as compared with known methods to detect asymmetry in the data distributions. 
Finally, this paper is concluded with a discussion in Section~\ref{sc5}. 
\section{Characterization  and Sufficiency}\label{sc2}
	In order to gain some insight into the purpose of this section,  the following illustrative example is offered.
	\subsection{One Sample Test of the Mean}\label{Ex1}
	Let $X_1, X_2$ be a random sample from a normal population with mean $\mu$ and variance $\sigma^2=1$. We consider testing $H_0:\,\mu=0$  versus $H_1:\,\mu=\delta$,  where $\delta$ is a fixed value. We present this simple toy example to illustrate our results shown in Section~\ref{sc22}, not to offer a contender to the usual t-test. In this example, the  statistic $\bar{X}=(X_1+X_2)/2$ can be used for MP testing. However, one may feel that, for example, the statistic $A=w X_1+(1-w)X_2$ could be reasonable for assessing the competing hypotheses $H_0$ and $H_1$, for some $w\in[0,1]$, $w\neq 0.5$. 
	 The vector $[\bar{X},A]^\top$ has a bivariate normal density function, with $\mathrm{E}\left(\bar{X}\right)=\mathrm{E}\left(A\right)=\mu$, $\mathrm{var}\left(\bar{X}\right)=\mathrm{cov}\left(\bar{X},A\right)=1/2$, and $\mathrm{var}(A)=w^2+(1-w)^2$. Then, by defining a joint density function of $\left(X_1,X_2\right)$-based  statistics $A_1, A_2$ in the form $f_{\mu}^{A_1,A_2}(u,v)$, 
	 	it is easy to observe that the ratio 
	$f_{\mu=\delta}^{\bar{X},A}(u,v)/f_{\mu=0}^{\bar{X},A}(u,v)$ does not depend on $v$, where $v$ is an argument of the joint density $f_{\mu}^{\bar{X},A}$ relating to $A$'s component. 	
	In particular, this means that after surveying the two data points $\bar{X}$ and $A$, we can improve the $\left(\bar{X},A\right)$-based decision-making mechanism by creating  the MP statistic 
	for testing $H_0$ vs. $H_1$  in the likelihood ratio form $f_{\mu=\delta}^{\bar{X},A}(\bar{X},A)/f_{\mu=0}^{\bar{X},A}(\bar{X},A)$,
	which only requires the computation of $\bar{X}$. 
	Thus, one might pose the question: can the observation above be generalized to extend Proposition~\ref{pr1}? In this case, it seems to be reasonable that to provide an essential property of the MP concept, relationships with other $D$-based statistics should be taken into account.

	The second aspect of our approach is to characterize a scenario where, say, statistic $A_1$ is more preferable 
	 in the construction of a test than  statistic $A_2$. In this context, as will be seen later, $A_1$ is superior to $A_2$, if  the ratio  $f_{\mu=\delta}^{A_1,A_2}(u,v)/f_{\mu=0}^{A_1,A_2}(u,v)$ does not depend on $v$. To exemplify the benefits of this rule, let us pretend that it is unknown that $\bar{X}$ is the best statistic in this subsection such that we can consider the following task. The problem then is to indicate a value of $a\in[0,1]$ in the statistic $T=a X_1+(1-a)X_2$ such that $T$ outperforms $A=w X_1+(1-w)X_2$, for all $w\in[0,1]$. Since the density  $f_{\mu}^{T,A}(u,v)$ is  bivariate normal,  simple algebra shows that  $f_{\mu=\delta}^{A_1,A_2}(u,v)/f_{\mu=0}^{A_1,A_2}(u,v)$ is not a function of $v$, 
	 if $a^2+(1-a)^2=aw+(1-a)(1-w)$. Then, the solution is  $a=0.5$. In reality, we do know that $T$, with $a=0.5$, is the best statistic in this framework. This example illustrates  our point.   

\subsection{Theoretical Results}\label{sc22}
In this section, the main results are provided in  Propositions~\ref{pr2}$-$\ref{pr5} below that establish the characterization of  MP tests. The proofs of Propositions~\ref{pr2}$-$\ref{pr5} are
included in the online supplementary materials for completeness and contain comments that augment the description of the obtained results. Proposition~\ref{pr5} revisits the characterization of  MP tests in the light of sufficiency.

To extend Proposition~\ref{pr1}, we define a joint density function of statistics $A=A(D)$ and $B=B(D)$ in the form $f^{A,B}_{k}(u,v)$, provided that $H_k$ is true, $k\in\{0,1\}$. 
Then, the likelihood ratio test statistic, $\Lambda$,  has the  following property.
\begin{proposition}\label{pr2}
	For any statistic $A=A(D)$, we have $f_{1 }^{\Lambda,A} (u,v)=u\,f_{0}^{\Lambda,A} (u,v)$, for all $u\ge 0$ and $v\in {{\mathbb{R}}^{1}}$.
\end{proposition}
Proposition~\ref{pr2}  emerges as a generalization of Proposition~\ref{pr1}, since 
$f_{1 }^{\Lambda,A} (u,v)=u{\rm \; }f_{0}^{\Lambda,A} (u,v)$ yields 
\begin{eqnarray*}
	f_{1 }^{\Lambda} (u)=\int f_{1 }^{\Lambda,A} (u,v) dv=u \,\int f_{0}^{\Lambda,A} (u,v) dv=u\,f_{0}^{\Lambda} (u).
\end{eqnarray*}

In the following claim, it is shown that Proposition~\ref{pr2} can be augmented  to imply a necessary and sufficient condition on test statistics distributions to present  MP decision-making techniques. Define $C$ to be a test threshold. 
\begin{proposition}\label{pr3}
	Assume a statistic for testing, $T=T(D)$, satisfies ${\Pr}_{k}(T\ge0)=1$, $k\in\{0,1\}$, one rejects $H_0$ when $T>C$.  The test statistic $T$ is MP if and only if  (iff) 
	$f_{1 }^{T,A} (u,v)=u\,f_{0}^{T,A} (u,v)$, for any statistic $A=A(D)$ and all  $u\ge 0$,  $v\in {{\mathbb{R}}^{1}}$.
\end{proposition}
Note that the condition "$T$ is strictly non-negative" is employed in~\cite{Thomas}, where a  monotonic   
logarithmic transformation of $T$, a test statistic, may improve the power of the $T$-based test when the corresponding TIE rate is asymptotically controlled at $\alpha$.
 However, the requirement 
$T\ge 0$ is not critical, because if we  evaluate a test statistic, say $G=G(D)$, that can be negative, then a monotonic transformation $T=g(G)\ge 0$ can assist in this case.

We can remark that, in scenarios where density functions of test statistics do not exist, the arguments employed in the proof of  Proposition~\ref{pr3} can be applied to obtain the next statement.
\begin{proposition}\label{pr33}
	The test statistic $T(D)>0$  is MP iff 
	$\mathrm{E}_{1 }\left\{g(D)\right\}=\mathrm{E}_{0 }\left\{g(D)T(D)\right\}$, for every function $g\in [0,1]$ of $D$. 	
\end{proposition}
\begin{remark} {\em 
Since $\mathrm{E}_{1 }\left\{g(D)\right\}=\mathrm{E}_{0 }\left\{g(D)\Lambda(D)\right\}$,
the condition $\mathrm{E}_{1 }\left\{g(D)\right\}=\mathrm{E}_{0 }\left\{g(D)T(D)\right\}$ implies 
	$\mathrm{E}_{0}\left[\left\{T(D)-\Lambda(D)\right\}g(D)\right]=0$, for every $g\in [0,1]$,
	which means, with probability one under $f_0$, we have $T=\Lambda$. Note also that, in Proposition~\ref{pr33}, we can use $g(D)$, satisfying $\mathrm{E}\left\{g(D)\right\}^m=\mathrm{E}\left\{g(D)\right\}$, for all $m>0$.   
	 }
\end{remark}	

The scheme used in the proof of Step (2)  of  Proposition~\ref{pr3} yields the following result.

\begin{proposition}\label{pr4}
	A statistic $T_{1}\ge 0 $ is more powerful than a statistic  $T_{2} $, if the ratio
	$f_{1}^{T_{1} ,T_{2}} (u,v)/f_{0}^{T_{1} ,T_{2}} (u,v)$ $= u$, for all  $u\ge 0$,  $v\in {{\mathbb{R}}^{1}}$.  
\end{proposition}
It is interesting to note that, by virtue of Propositions~\ref{pr1} and~\ref{pr3},  for any $D$-based statistic $A=A(D)$, we have  $f_{1 }^{T,A} (u,v)={f_{0}^{T,A} (u,v)}\, u$ and then 

\noindent $f_{1 }^{T,A} (u,v)$ $={f_{0}^{T,A} (u,v)} f_{1 }^{T} (u)/f_{0}^{T} (u),$ 
if $T$ is MP. That is to say, $f_{1 }^{A|T} (u,v)={\rm \; }f_{0}^{A|T} (u,v)$, where the notation $f_{k}^{A|T}$ means a conditional density function of $A$ given $T$ under $H_k$, $k\in\{0,1\}$. In this case, when $A$ is independent of $T$, we obtain  $f_{1 }^{A} (v)=f_{0}^{A} (v)$, and then $A=A(D)$ cannot discriminate the hypotheses.
We can write that  $A=A(D)$ is ancillary, meaning $f_{1 }^{A} =f_{0}^{A}$. 
This motivates us to associate the results above with the principle of sufficiency.   

According to~\cite{Bah}, in the considered framework, we can call $T=T(D)$ to be a sufficient test statistic, if $f_{1 }^{A|T} (u,v)={\rm \; }f_{0}^{A|T} (u,v)$,  for each $A=A(D)$ and all  $u\ge 0$,  $v\in {{\mathbb{R}}^{1}}$. In this context, the statements mentioned above assert  the next result.
\begin{proposition}\label{pr5} 
	The following claims are equivalent:
	\begin{enumerate}
		\item[{\rm (i)\phantom{i}}]
		$T=T(D)$ is sufficient and $f_{1 }^{T} (u)/f_{0}^{T} (u)=u$;
		\item[{\rm (ii)\phantom{i}}]		
		$T$ is a MP statistic for testing the competing hypotheses $H_0$ and $H_1$.
	\end{enumerate}
\end{proposition}

	 Proposition~\ref{pr5} presents an argument to the reasonableness of making a statistical inference  based solely on the corresponding sufficient statistics.  
\begin{remark}\label{rem2}
{\em	We can note the following facts: 
	\begin{enumerate}
		\item[{\rm (i)\phantom{i}}]
		\cite{Kagan} have exemplified a sufficiency paradox, when an insufficient statistic preserves the Fisher information. 
		\item[{\rm (ii)\phantom{i}}]
		In order to  extend Proposition~\ref{pr5}, statements related to a wide spectrum of  Basu's theorem-type results (e.g.,~\citealp{Basu}) can be employed, in certain situations. 
		To the best of  our knowledge, there are no direct applications of Basu’s theorem to the questions considered in the present article.
		\item[{\rm (iii)\phantom{i}}] In Bayesian styles of  testing   (e.g.,~\citealp{Bayes}), Proposition~\ref{pr1} can be extended to treat Bayes Factors, see, e.g., Proposition 5 of \cite{VexlerP}, in this context. Then, Propositions~\ref{pr2} and~\ref{pr3} can be easily modified to establish integrated MP tests with respect to incorporated prior information~\citep{VexlerM}.		
	\end{enumerate}
}
\end{remark}


\section{Applications}\label{sc3} Section~\ref{sc2}
carries out the relatively  general underlying theoretical framework for the MP characterization concept. In this section, we outline three applications of the proposed MP characterization principle, by modifying well-accepted statistical tests in an easy to implement manner.
It is hoped that  the proposed MP characterization can provide different benefits for developing, improving, and comparing decision-making algorithms in statistical practice.   

 A common problem arising in statistical inference is the need for methods to modify a given test statistic in order to improve the performance of controlling the TIE rate  
and  power of the corresponding decision-making scheme. For example, 
 the accuracy of asymptotic approximations for the null distribution of a test statistic 
 may be increased by incorporating Bartlett correction type mechanisms or/and location adjustment techniques. In this context, we refer the reader to the following examples: \cite{hall1990methodology},    for modifying nonparametric empirical likelihood ratios; \cite{Chen}, for different transformations of the t-test statistics assessing the mean of asymmetrical distributions. Recently,
 \cite{Thomas} proposed to  use a logarithmic transformation  to obtain a potential  increase in power of the transformed statistic-based test.

This section  demonstrates use of the considered MP principle, following the simple idea outlined below. Suppose that we have a reasonable test statistic $T_o$ and we wish to improve \(T_o\) to be in a form, say $T_N$,  approximately satisfying the claim $f_{1 }^{T_N,A} (u,v)=u\,f_{0}^{T_N,A} (u,v)$, for any statistic $A=A(D)$ and all  $u\ge 0$,  $v\in {{\mathbb{R}}^{1}}$. Given that in general 
nonparametric settings there are no MP tests, it would be attractive to reach the   
MP property  $f_{1 }^{T_N,A} (u,v)=u\,f_{0}^{T_N,A} (u,v)$ at least for some statistic $A$, especially for some ancillary statistic. Informally speaking, by having $A$ with $f_{1}^{A} = f_{0}^{A}$, we can remove the influence of $A$ from $T_o$ to create $T_N$ such that 
  the ratio $f_{1 }^{T_N,A} (u,v)/f_{0}^{T_N,A}(u,v)$ is a function of $u$ only.
 In this case, Proposition~\ref{pr4} could insure that $T_N$ outperforms $A$.  
This can be achieved via an independence between $T_N$ and $A$  that is exemplified in 
Sections~\ref{sc31}, \ref{sc32}, and~\ref{sc33}  in detail.

	Through the following examples, we  aim to show our approach in   an   intuitive manner. 
\subsection{Examples of the Use of Ancillary Statistics}\label{sc311}
We begin with displaying the  toy examples below that illustrate our key
idea. 

\begin{figure}
	\begin{center}
		\includegraphics[width=5in]{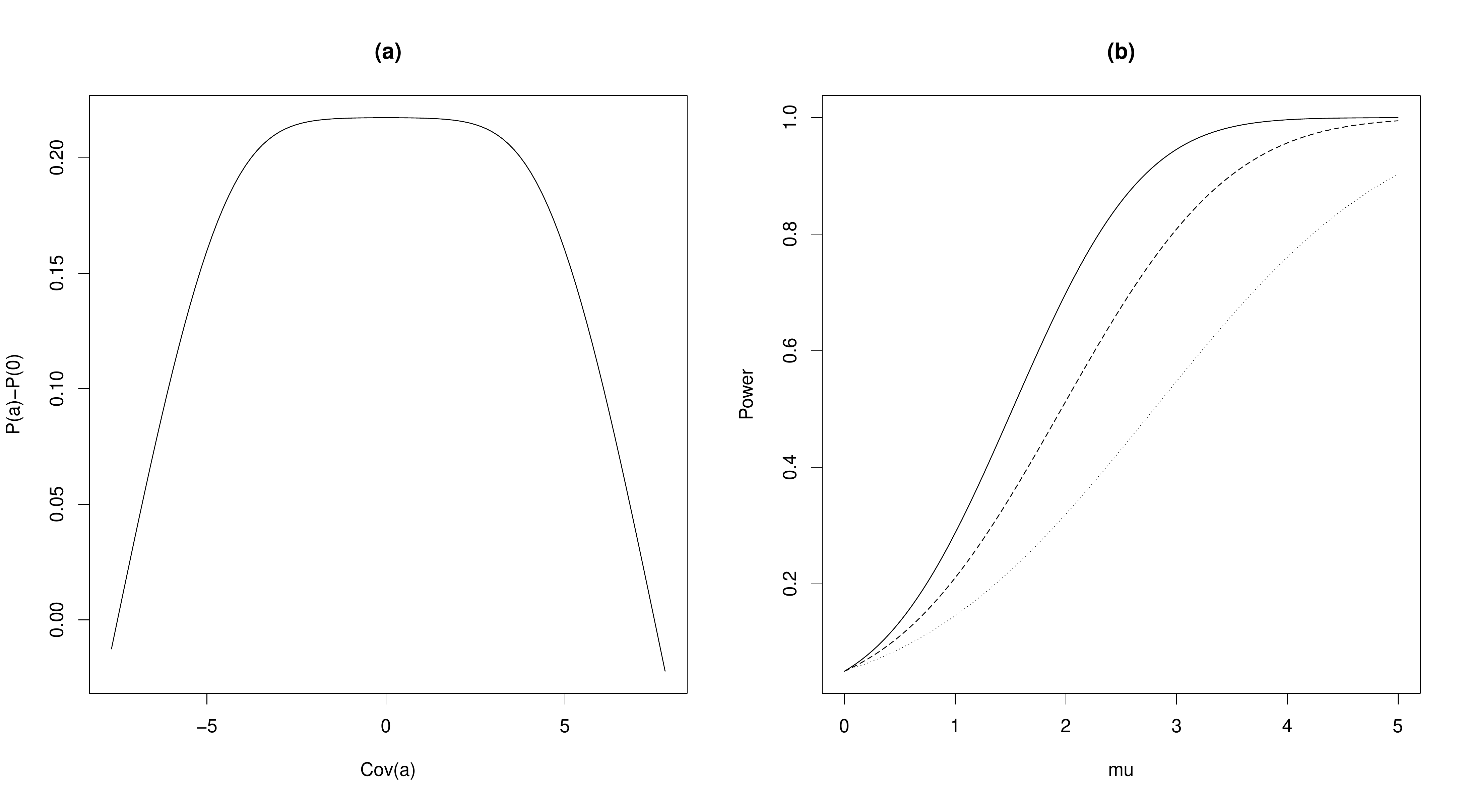}
		\caption{Graphical evaluations related to the examples shown in Section~\ref{sc311}. Panel (a) plots $P(a)-P(0)$, the power of the $\left(T+a(X_1-X_2)\right)$-based test minus the power of the $T$-based test  at the $\alpha=0.05$ level, against the covariance $Cov(a)=\mathrm{cov}\left(T+a(X_1-X_2),X_1-X_2\right)$, for $a\in[-0.01,0.9]$, where $T=0.5(X_1+X_2)$, $X_1\sim N(\mu,1)$, $X_2\sim N(\mu,4^2)$, $\mathrm{E}_0(X_i)=0$,  $\mathrm{E}_1(X_i)=5$, $i\in[1,2]$.  Panel (b) plots the powers
		$P_{T_O}(\mu)={\Pr}_{\mu}\left\{T_O>C^{T_O}_{\alpha}\right\}$ (solid line),		
		$P_{T_N}(\mu)={\Pr}_{\mu}\left\{T_N>C^{T_N}_{\alpha}\right\}$ (longdashed line),
		and	
		$P_{T}(\mu)={\Pr}_{\mu}\left\{T>C^T_{\alpha}\right\}$ (dotted line) at the $\alpha=0.05$ level,
        where 
        $T=(X_1+X_2+X_3)/3$, $T_N=T+\gamma(X_1-X_2)$, $T_O=\left(X_1/\sigma^2_1+X_2/\sigma^2_2+X_3/\sigma^2_3\right)/\left(1/\sigma^2_1+1/\sigma^2_2+1/\sigma^2_3\right)$, $\gamma=\left(\sigma^2_2-\sigma^2_1\right)\left(\sigma^2_2+\sigma^2_1\right)^{-1}/3$, $X_1\sim N(\mu,1)$, $X_2\sim N(\mu,4^2)$, and $X_3\sim N(\mu,3^2)$, for $\mu\in [0,5]$.
    }
		\label{Example}
	\end{center}
\end{figure}

Let independent data points $X_1$ and $X_2$ be observed; when it is assumed that $X_i\sim N\left(\mu,\sigma^2_i\right)$, $i\in [1,2]$ and $\sigma^2_1\neq \sigma^2_2$ are known.  One can use the simple  statistic $T=0.5(X_1+X_2)$ to test  $H_0:$ $\mu=0$ against $H_1:$ $\mu>0$.  Easily, one can confirm  that $X_1-X_2$ is an ancillary statistic. We now consider a mechanism for transforming $T$ and making a modified test statistic that is independent of $X_1-X_2$.
Define $T_N=T+\gamma(X_1-X_2)$,  a transformed version of $T$, where $\gamma$ is a root of the equation $\mathrm{cov}\left(T_N,X_1-X_2\right)=0$. Then, we obtain $\gamma=0.5\left(\sigma^2_2-\sigma^2_1\right)/\left(\sigma^2_2+\sigma^2_1\right)$. Thus, the derived statistic 
\begin{eqnarray*}
T_N&=&T+\frac{\sigma^2_2-\sigma^2_1}{2\left(\sigma^2_2+\sigma^2_1\right)}\left(X_1-X_2\right)
=\frac{\sigma^2_2X_1+\sigma^2_1X_2}{\sigma^2_2+\sigma^2_1}
\,\,=\,\,\frac{X_1/\sigma^2_1+X_2/\sigma^2_2}{1/\sigma^2_1+1/\sigma^2_2}
\end{eqnarray*}
is certainly a successful transformation of the initial statistic $T$, which presents the MP test statistic. 
For instance, we denote the power 
\begin{eqnarray*}
	P(a)={\Pr}_{\mu=5}\left\{T+a(X_1-X_2)>C(a)\right\},\, C(a):\, {\Pr}_{\mu=0}\left\{T+a(X_1-X_2)>C(a)\right\}=0.05, 
\end{eqnarray*}
when $\sigma_1=1$ and $\sigma_2=4$.
Figure~\ref{Example}(a) depicts the function $P(a)-P(0)$, the difference between the power levels of the 
$\left(T+a(X_1-X_2)\right)$-based test and those of the $T$-based test at $\alpha=0.05$, plotted against the function $Cov(a)=\mathrm{cov}\left(T+a(X_1-X_2),X_1-X_2\right)$, for $a\in[-0.01,0.9]$. 
As expected,  the function $P(a)-P(0)$ reaches its maximum when $Cov(a)=0$. Moreover, it turns out that we do not need much accuracy in approximating the equation $Cov(a)=0$ to outperform the $T$-based test when we use the modified test statistic $T+a(X_1-X_2)$. 
Then,  intuitively, we can suppose that
a transformed test statistic could include estimated elements while still providing good power characteristics for its decision-making algorithm.

	In various situations, we shall  not exclude the possibility that there exists more than one
	ancillary statistic for a given testing statement.  Let us exemplify such case, assuming we observe
	$X_1$, $X_2$, and $X_3$ from the normal distributions  $N\left(\mu,\sigma^2_1\right)$, $N\left(\mu,\sigma^2_2\right)$,  and $N\left(\mu,\sigma^2_3\right)$, respectively, where $\sigma^2_i$, $i\in [1,2,3]$, are known. Suppose we are interested in testing $H_0:$ $\mu=0$ vs. $H_1:$ $\mu>0$. The statistic to be modified is $T=(X_1+X_2+X_3)/3$. The observation $X_1-X_2$ is an ancillary statistic with respect to $\mu$. Define $T_N=T+\gamma(X_1-X_2)$ with $\gamma=\left(\sigma^2_2-\sigma^2_1\right)\left(\sigma^2_2+\sigma^2_1\right)^{-1}/3$, thereby obtaining that
	$\mathrm{cov}\left(T+\gamma(X_1-X_2),X_1-X_2\right)=0$. Then, it is clear that $T_N$ is somewhat better than $T$, but $T_O=\sum_{i=1}^3\left(X_i/\sigma^2_i\right)/\sum_{i=1}^2\left(1/\sigma^2_i\right)$ is superior to $T_N$ in the terms of this example. Define the powers $P_{T}(\mu)={\Pr}_{\mu}\left\{T>C^T_{0.05}\right\}$, $P_{T_N}(\mu)={\Pr}_{\mu}\left\{T_N>C^{T_N}_{0.05}\right\}$, and $P_{T_O}(\mu)={\Pr}_{\mu}\left\{T_O>C^{T_O}_{0.05}\right\}$, where the test thresholds $C^{T}_{0.05}$, $C^{T_N}_{0.05}$, and $C^{T_O}_{0.05}$ satisfy  	
	$
		{\Pr}_{\mu=0}\left\{T>C^T_{0.05}\right\}= {\Pr}_{\mu=0}\left\{T_N>C^{T_N}_{0.05}\right\}=
		{\Pr}_{\mu=0}\left\{T_O>C^{T_O}_{0.05}\right\}=
		0.05. 
	$ Figure~\ref{Example}(b) exemplifies the behavior of the functions  $P_{T}(\mu)$, $P_{T_N}(\mu)$, and $P_{T_O}(\mu)$, when  $\sigma_1=1$, $\sigma_2=4$, and $\sigma_3=3$.

	Thus, to improve a given test statistic, say $T$, we can suggest that one pays  attention to a relevant ancillary statistic, say $A$,  modifying $T$ to be independent (or approximately independent) of $A$.

	Note that, although the concept of ancillarity asserts that ancillary statistics do not provide information about the parameters of  interest, different roles of ancillary statistics in parametric estimation  have been dealt with extensively in the literature. In this context, for an extensive review, we refer the reader to~\cite{Gh2}. 
	For example, assume we observe the vectors  $[X_i,Y_i]^\top$, $i\in[1,\ldots,n]$, from a bivariate normal distribution with $\mathrm{E}(X_1)=\mathrm{E}(Y_1)=0$, 
	$\mathrm{var}(X_1)=\mathrm{var}(Y_1)=1$, and $\mathrm{corr}(X_1,Y_1)=\rho$,	
	where $\rho\in (-1,1)$ is unknown. The statistics $U_1=\sum_{i=1}^n X^2_i$ and $U_2=\sum_{i=1}^n Y^2_i$ are ancillary. According to ~\cite{Gh2}, to define unbiased estimators of $\rho$, it can be recommended  to use the statistics $\sum_{i=1}^{n}X_iY_i/U_j$, $j\in[1,1]$. As another example, when ancillary statistics are applied, we outline a case of so-called Monte Carlo swindles,   simulation based methods	
	 that allow  small numbers of generated samples to produce statistical accuracy at the level one would expect from  much larger numbers of generated samples. \cite{boos1998applications} discussed the following procedure. To estimate the variance of the sample median $M$ of a normally distributed sample $X_1,\ldots,X_n$, the Monte Carlo swindle approach estimates $\mathrm{var}\left(M -\bar{X}\right)$ (instead of $\mathrm{var}\left(M\right)$) by using the $N$ Monte Carlo samples of $X_1,\ldots,X_n$
	  	  and then  $\mathrm{var}\left(\bar{X}\right)=\mathrm{var}\left(X_1\right)/n$ is added to obtain an efficient  estimate of $\mathrm{var}\left(M\right)$. 
	  	  In order to justify this framework, we employ that the statistic $V = (X_1 - \bar{X},\ldots,X_n - \bar{X})$ is ancillary. Now, since $\bar{X}$ is complete sufficient, $\bar{X}$ and $V$ are independent by Basu's theorem. Then, $\bar{X}$ is independent of the sample median of $V$. Therefore,
	  	  \begin{eqnarray*}
	  	  	\mathrm{var}\left(M\right)\,=\,	\mathrm{var}\left(M-\bar{X}+\bar{X}\right)\,=\,
	  	  	\mathrm{var}\left(M-\bar{X}\right)+\mathrm{var}\left(\bar{X}\right).	  	  	
	  	  \end{eqnarray*}	  	  	  	    
	  	   It is clear that, when $X_1$ has a normal distribution,  the contribution from $\mathrm{var}\left(X_1\right)/n$  to $\mathrm{var}\left(M\right)$ is much larger than the contribution from $\mathrm{var}\left(M-\bar{X}\right)$, where the component 
	  	   $\mathrm{var}\left(M-\bar{X}\right)$ is proposed to be  estimated by simulation. 
	  	   This limits the error in estimation by simulation  to a  small part of $\mathrm{var}\left(M\right)$ (for  details, see \citealp{boos1998applications}).

\subsection{Theoretical Support}
The  point of view mentioned above can be supported by the following results. 
Assume we have a test statistic $Y=Y(D)$, and the ratio $L^{Y}(u)=f^{Y}_1(u)/f^{Y}_0(u)$ is a monotonically increasing function that has an inverse function, say $W(u)$. In this scenario, $Y$ can be transformed into the form $Y_N=L^{Y}(Y)$, thereby implying that
 \begin{eqnarray*}
 	f^{Y_N}_k(u)\,=\,\frac{d}{du}{\Pr}_k\left\{L^{Y}(Y)\le u\right\}\,=\,\frac{d}{du}{\Pr}_k\left\{Y\le W(u)\right\}\,=\,f_k^{Y}\left(W(u)\right)\frac{d}{du}W(u),\, k\in\{0,1\}.
 \end{eqnarray*}
This means that the likelihood ratio 
 \begin{eqnarray*}
	f^{Y_N}_1\left(Y_N\right)/f^{Y_N}_0\left(Y_N\right)\,=\,f_1^{Y}\left(W(Y_N)\right)/f_0^{Y}\left(W(Y_N)\right)\,=\,f_1^{Y}\left(Y\right)/f_0^{Y}\left(Y\right)\,=\,Y_N.
\end{eqnarray*}
(See  Proposition~\ref{pr1}, in this context.) Then, we  state the next proposition.
Let a statistic $A=A(D)$ satisfy  $f_{1}^{A} = f_{0}^{A}$. Suppose we have the decision-making procedure based on a statistic $T=T(D)$, and we can modify $T$ to  $T_N$ to achieve  $T_N$ and $A$ as independent terms under $H_0$ and $H_1$, when $T=\psi\left(T_N,A\right)$, for a bivariate function $\psi$. We conclude that:
\begin{proposition}\label{pr6} 
	The $T_N$-based test considered above is superior to that based on $T$, if the ratio  
	$L^{T_N}(u)=f^{T_N}_1(u)/f^{T_N}_0(u)$ is a monotone  function.		
\end{proposition}
The proof is deferred to the online supplementary materials.

	Note that Proposition~\ref{pr6} gives some insight into  the connection between the power of statistical tests and ancillarity,  concepts that seem to be unrelated, since ancillary statistics  cannot solely discriminate the competing hypotheses $H_0$ and $H_1$.  

Proposition~\ref{pr6} depicts the rationale for modifying the following well-known test statistics.
   
\subsection{One Sample t-Test for the Mean} \label{sc31}
Assume we observe independent and identically distributed (i.i.d.) data points $X_1,X_2,\ldots,X_n$ that provide $D=\left\{X_1,\ldots,X_n\right\}$. For testing the hypothesis $H_0:$ $\mu=0$ vs. $H_1:$ $\mu>0$, where $\mu=\text{E}X_1$, the well-accepted statistic is $T_o=n^{0.5}\bar{X}/\sigma$, where $\bar{X}=\sum_{i=1}^nX_i/n$, and $\sigma^2=\text{var}(X_1)$. 

In this testing statement, it seems that the statistic $S^2=\sum_{i=1}^n\left(X_i-\bar{X}\right)^2/(n-1)$, the sample variance, is approximately ancillary with respect to $\mu$. 
Note also that  $T_o$ and $S^2$ are independent, when  $X_1\sim N\left(\mu,\sigma^2\right)$ ($T_o$ is MP, in this case). Then, we denote the statistic $T_N(\gamma)=T_o+\gamma S^2$ and derive a value of $\gamma$, say $\gamma_0$, that insures $\text{cov}\left(T_N(\gamma_0),S^2\right)=0$. 
The statistic $T_N(\gamma_0)$ is a basic ingredient of the modified test statistic we will propose.
To this end, we define $\mu_k=\text{E}\left(X_1-\mu\right)^k,\, k=3,4$, and employ the results from \cite{O'Neill} in order to obtain 
\begin{eqnarray*}
	\gamma_0&=&\sigma^{-1}\mu_3n^{0.5}/\left(\sigma^4\frac{n-3}{n-1}-\mu_4\right)\,\approx \,
	-\sigma^{-1}\mu_3n^{0.5}/\text{var}\left\{\left(X_1-\mu\right)^2\right\}.
\end{eqnarray*} 

The additional argument for using $T_N(\gamma_0)$ in a test for $H_0$ vs. $H_1$ can be explained in the following simple fashion.
 It is clear that the stated testing problem can be treated in the context of  a confidence interval estimation of $\mu$. Thus, there is a relationship between the quality of testing  $H_0:$ $\mu=0$ and the variance of an estimator of $\mu$ involved in corresponding decision-making schemes. (For example, the t-test, $T_o$, uses $\bar{X}$ to estimate $\mu$.)   For the sake of simplicity, consider 
			$
			\tilde{T}_N(\gamma)=\bar{X}+\gamma\left(S^2-\sigma^2\right)
			$
			that satisfies $\mathrm{E}\tilde{T}_N(\gamma)=\mu$ and $\mathrm{var}\left\{\tilde{T}_N(0)\right\}$ $=\mathrm{var}\left(\bar{X}\right)$. 
			To find a value of $\gamma$ that minimizes $\mathrm{var}\left\{\tilde{T}_N(\gamma)\right\}$, we can solve the  equation $d\left[\mathrm{var}\left\{\tilde{T}_N(\gamma)\right\}\right]/d\gamma=0$, where it is assumed we can write 
   \begin{eqnarray*}
    \frac{d}{d\gamma}\mathrm{var}\left\{\tilde{T}_N(\gamma)\right\}&=&  \frac{d}{d\gamma}\mathrm{E}\left\{\bar{X}+\gamma \left(S^2-\sigma^2\right)-\mu\right\}^2 =
    \mathrm{E} \frac{d}{d\gamma}\left\{\bar{X}+\gamma \left(S^2-\sigma^2\right)-\mu\right\}^2
    \\
    &=&2\mathrm{E}\left\{\bar{X}+\gamma \left(S^2-\sigma^2\right)-\mu\right\}\left(S^2-\sigma^2\right).
   \end{eqnarray*} 
   Then, the root
   \[
			\gamma=-\mathrm{cov}\left(\bar{X},S^2\right)/\mathrm{var}\left(S^2\right) 
			=\mu_3/\left\{\sigma^4\left(n-3\right)/\left(n-1\right)-\mu_4\right\} 
			\approx 
			-\mu_3/\mathrm{var}\left\{\left(X_1-\mu\right)^2\right\}
	\]
 minimizes $\mathrm{var}\left\{\tilde{T}_N(\gamma)\right\}$  and is  $\gamma_0\sigma/n^{0.5}$, where $\gamma_0$ implies $\mathrm{cov}\left(T_N(\gamma_0),S^2\right)=0$. That is, 
			$ \mathrm{var}\left\{\tilde{T}_N\left(\gamma_0\sigma/n^{0.5}\right)\right\}$ $\le$ $\mathrm{var}\left(\bar{X}\right)$.   
			The statistic $T_N(\gamma_0)$ includes $\bar{X}$ multiplied by $n^{0.5}/\sigma$.  This confirms that the statistic $T_N(\gamma_0)$ can  be  somewhat more powerful than $T_o$, 
			in the terms of testing $H_0$ vs. $H_1$, for various scenarios of $X_1$'s distributions.

Finally, we standardize the test statistic $T_N(\gamma_0)$  to be able to control its TIE rate, denoting the $\alpha$ level  decision-making rule:  the null hypothesis $H_0$ is rejected if 
\[T_N=\hat{\Delta}^{-0.5}\left[T_o-\frac{\hat{\mu}_3n^{0.5}}{\sigma\,\hat{\text{var}}\left\{\left(X_1-\mu\right)^2\right\}}\left(S^2-\sigma^2\right)\right]>z_\alpha,\]
where $\hat{\mu}_3=n^{-1}\sum_{i=1}^n\left(X_i-\bar{X}\right)^3$ and $\hat{\text{var}}\left\{\left(X_1-\mu\right)^2\right\}$ 
$=n^{-1}\sum_{i=1}^n\left\{\left(X_i-\bar{X}\right)^2-\sigma^4\right\}^2$ estimate  $\mu_3$ and
$\text{var}\left\{\left(X_1-\mu\right)^2\right\}$, respectively; 
$$\hat{\Delta}=1-S^{-2}\hat{\mu}_3^2\left[n^{-1}\sum_{i=1}^n\left\{\left(X_i-\bar{X}\right)^2-S^4\right\}^2\right]^{-1}$$
is the sample estimator of $\Delta=\text{var}\left(T_N(\gamma_0)\right)$; and the threshold $z_\alpha$ satisfies $\Pr(Z>z_{\alpha})=\alpha$ with $Z\sim N(0,1)$. It can be interesting to rewrite $T_N$ in the form $T_N=n^{0.5}\bar{Y}/\sigma$, where $Y_i=\hat{\Delta}^{-0.5}\left[X_i-\hat{\mu}_3\left\{\left(X_i-\bar{X}\right)^2n/(n-1)-\sigma^2\right\}/\hat{\text{var}}\left\{\left(X_1-\mu\right)^2\right\}\right]$. The statistic $T_N$ is asymptotically $N(0,1)$-distributed, under $H_0$. Certainly, in the case of $X_1\sim N(\mu,\sigma^2)$, meaning  $T_o$ is MP, we have  
$\mu_3=0$. In the form  $T_N$, an adjustment for the skewness of the underlying data is in effect.    

	Note that, in the transformation of the test statistic shown above, we achieve uncorrelatedness between $T_N$ and $S^2$,  thus simplifying the development of the nonparametric procedure.  It is clear that the equality $\text{cov}\left(T_N,S^2\right)=0$ is  essential to the asymptotic independence between $T_N$ and $S^2$ (e.g., \citealp[pp. 181–206]{ghosh2021advances}). 

Section~\ref{sc4} uses  extensive Monte Carlo evaluations  to demonstrate an efficiency of the statistic $T_N$ for testing $H_0$ against various alternatives.

The testing procedures based on $T_o$ and $T_N$ require $\sigma$ to be known. 
This restriction can be overcome, for example by using bootstrap type strategies, Bayesian techniques and/or p-value-based methods introduced by~\cite{Pvalue}. In this article, we only note that there are practical applications in which it is reasonable to assume $\sigma$ is known, e.g., \cite{Data}, \citet[Section 3.2]{boos1998applications},  as well as  \citet[p. 1729]{Bayes}. In many biostatistical studies, biomarkers values are scaled  in such a way that their variance $\sigma^2=1$. 
We also remark  that developments of simple test statistics improving the t-test, $T_o$, can be of a theoretical  interest.  
\subsection{One Sample Test for the Median}\label{sc32}
A sub-problem related to comparisons between mean and quantile effects can be considered as follows:
Let $X_1,X_2\ldots,X_n$ be continuous i.i.d. observations with $\mathrm{E}X_1=0$. We are interested in testing the hypothesis
$H_0:$ $\nu=0$ vs. $H_1:$ $\nu>0$, where $\nu$ denotes the median of $X_1$.
 This statement of the problem can  be  found in various practical applications related to testing linear regression residuals as being symmetric, and  pre-and post-placebo  paired comparison of biomarker measurements
  as well as, for example, when researchers investigate data associated with  radioactivity detection in drinking water, where the population mean is known; see Section 4 in~\cite{Data2}.

To test for  $H_0$, it is reasonable to use a statistic in the form 
$$T_o=2n^{0.5}X_{(n/2)} f\left(X_{(n/2)}\right),$$ 
where $X_{(n/2)}$ is the sample estimator of $\nu$ based on the order  statistics $X_{(1)}<X_{(2)}< \cdots< X_{(n)}$, and $1/f\left(X_{(n/2)}\right)$ is a measure of scale, with 
$f(u)$ being the density function of $X_1$. The statistic $A=n^{0.5}\bar{X}/\sigma$ can be selected as  approximately ancillary with respect to $\nu$, since $\mathrm{E}X_1=0$.  The known facts we use are: (a) the statistics $T_o$ and $A$ have  an asymptotic bivariate normal distribution with the parameters shown in~\cite{Median}; and 
(b) if $X_1\sim N(0,1)$ then Basu's theorem asserts that  $\bar{X}$ and $X_{(n/2)}-\bar{X}$ are independent, since $\bar{X}$ is a complete sufficient statistic and $X_{(n/2)}-\bar{X}$ is ancillary.     
That is to say, in a similar manner to the development shown in Section~\ref{sc31}, we can improve $T_o$ by focusing on the statistic 
\[
T=\left\{\sigma T_o /\mathrm{E}\left|X_1-\nu\right|-A\right\}\left\{\sigma^2\left(\mathrm{E}\left|X_1-\nu\right|\right)^{-2}-1\right\}^{-0.5}
\] 
that satisfies $\mathrm{var}\left(T\right)\to 1$ and $\mathrm{cov}\left(T,A\right)\to 0$, as $n\to \infty$. Thus, the test we propose is as follows: to reject $H_0$, if
\[T_N=\left\{2n^{0.5}X_{(n/2)} \hat{f}\left(X_{(n/2)}\right)\, S/\hat{w}-n^{0.5}\bar{X}/S\right\} \left\{S^2\hat{w}^{-2}-1\right\}^{-0.5}\,\, >\, \, z_\alpha,\]
where $S^2=\sum_{i=1}^n\left(X_i-\bar{X}\right)^2/(n-1)$, $\hat{w}=\sum_{i=1}^n\left|X_i-X_{(n/2)}\right|/n$, and $\hat{f}$ means a kernel estimator of the density function $f$. In order to estimate $f\left(X_{(n/2)}\right)$, 
we suggest employment of the R-command (\citealp{R}): $$\mathtt{density(X,from=median(X),to=median(X))[[2]][[1]]}$$

It is clear that, for two-sided testing $H_0:$ $\nu=0$ vs. $H_1:$ $\nu\ne 0$, we can apply the rejection rule: $T_N^2>\chi^2_1(\alpha)$, where the threshold $\chi^2_1(\alpha)$ satisfies $\Pr\left\{Z>\chi^2_1(\alpha)\right\}=\alpha$ with a random variable $Z$ having  a chi-square distribution with one degree of freedom.

It can be  remarked that the testing algorithm shown in this section can be easily extended to make decisions regarding quantiles of underlying data distributions; see Section~\ref{sc5}, for details.  
\subsection{Test for the Center of Symmetry}\label{sc33} 
In many practical applications, e.g., paired testing for pre-and post-treatment effects, we may be interested in testing that the center of symmetry of the paired observations is zero. 
To this end, we assume that i.i.d. observations $X_1,\ldots,X_n$ are from an unknown symmetric distribution $F$ with $\sigma^2=\text{var}\left(X_1\right)<\infty$. According to \cite{bickel2012descriptive}, the natural location
parameter, say $\nu$, for $F$ is its center of symmetry.  We are interested in testing the hypothesis  
$H_0:$ $\nu=0$ vs. $H_1:$ $\nu>0$.

The  statistics $T_o=n^{0.5}\bar{X}/S$ and $T_1=2n^{0.5}X_{(n/2)} \hat{f}\left(X_{(n/2)}\right)$  
 are reasonable to be employed for testing $H_0$, where $S^2=\sum_{i=1}^n\left(X_i-\bar{X}\right)^2/(n-1)$ and $\hat{f}(u)$ estimates $f(u)=dF(u)/du$. 
The components of $T_o$ and $T_1$ are specialized in Sections~\ref{sc31} and~\ref{sc32}. It is clear that if $F$ were known to be a normal distribution function, then $T_o$ outperforms $T_1$, whereas when $F$ were known to be  a distribution  of, e.g., the random variable $\xi_1-\xi_2$, where $\xi_1,\xi_2$ are independent and identically $Exp(1)$-distributed, $T_1$ outperforms $T_o$. 

Consider, for example, the statistic $T_o$ as a test statistic   to be modified and the statistic $A=\left(\bar{X}-X_{(n/2)}\right)D^{-0.5}\,$ having a role of an approximately ancillary statistic, where $D=\sigma^2-\text{E}|X_1-\nu|/f(\nu)+1/\left(4f^2(\nu)\right)$. Then, following the concept and the notations defined in Sections~\ref{sc31} and~\ref{sc32}, we  propose to reject $H_0$, if
\begin{eqnarray*}
	T_N&=& \left\{T_o+\delta n^{0.5}\left(\bar{X}-X_{(n/2)}\right)\hat{D}^{-0.5}\right\}V^{-0.5}>z_{\alpha},
	\\
	\hat{D}&=&S^2-\hat{w}/\hat{f}\left(X_{(n/2)}\right)+\left(4\hat{f}^2\left(X_{(n/2)}\right)\right)^{-1},\,\,
	\delta\,=\,\left\{\hat{w}/\left(2S\hat{f}\left(X_{(n/2)}\right)\right)-S\right\}\hat{D}^{-0.5},
	\\
	V&=&1+\delta^2+\frac{2\delta S}{\hat{D}}- \frac{\delta \hat{w}} {\hat{D}^{0.5}S\hat{f}\left(X_{(n/2)}\right)}
	\,=\,1+\frac{2\hat{w}}{\hat{D}\hat{f}\left(X_{(n/2)}\right)}
	-\frac{1}{\hat{D}}\left(\frac{\hat{w}}{2S\hat{f}\left(X_{(n/2)}\right)}+S\right)^2.
\end{eqnarray*}         
We will experimentally demonstrate that the $T_N$-based test can  combine attractive power properties of the $T_o$-and $T_1$-based tests.

\begin{remark}\label{rem4}
{\em Note that, in this section, the test statistics $T_N$ are targeted to improve the statistics $T_o$. In this section's framework, there are no MP decision-making mechanisms. Thus, in general, it can be assumed  we can find  decision-making procedures that outperform  the $T_N$-based tests in certain situations. }
\end{remark}

Section~\ref{sc4} numerically examines properties of the decision-making schemes derived in this section. 
\section{Numerical Simulations}\label{sc4}
We conducted a Monte Carlo study to explore the performance of the proposed  transformations of 
the tests about the mean, median, and center of symmetry as described in  Section~\ref{sc3}.
 In terms related to evaluations of nonparametric decision-making procedures, it can be noted  that  there are no MP tests, in the frameworks of Sections~\ref{sc31}, \ref{sc32}, and \ref{sc33}. 
	We therefore compare the tests based on the given statistics, $T_o$, with those based on the  corresponding statistics $T_N$, the modifications of $T_o$, under various designs of $H_0$/$H_1$-underlying data distributions.     
	The aim of the numerical study is to confirm that the proposed method can provide improvements in the context of statistical power. In Sections~\ref{sc42} and \ref{sc43},  for additional comparisons, we demonstrate 
	the Monte Carlo power of the one-sample Wilcoxon-Mann-Whitney test that is frequently  used in  applications, where researchers are interested in assessing  the hypothesis $H_0:$ $\nu=0$ when $\nu$ is the median of observations. Note that, in practice, it is very  difficult to find a nonparametric alternative to the  one sample  t-test for the mean. Then, in Section~\ref{sc41}, where the t-test and its transformation defined in Section~\ref{sc31} are evaluated, we include a bootstrapped (nonparametric resampling) version of the original t-statistic $T_o$ to be compared with the corresponding statistic $T_N$, expecting that the bootstrapped t-test may outperform the original t-test in several nonparametric scenarios ~\citep{efron1992bootstrap}.

To evaluate the tests, we generated $55,000$ independent samples of size \( n\in \{n_1,\ldots,n_J\} \)  from different distributions corresponding to, say, designs $D_{km}$, $k\in\{0,1\}$, $m\in\{1,\ldots,M\}$. In this scheme, designs $D_{km}$,  $m\in\{1,\ldots,M\}$, fit  hypotheses $H_k$, $k\in\{0,1\}$, respectively.
	Each of the presented bootstrap simulation
	results are based on  $55,000$ replications with $1,000$ bootstrap samples.

Let the notation $T\left(D_{km}\right)$ represent a test statistic $T$ conducted with respect to design $D_{km}$, $k\in\{0,1\}$, $m\in\{1,\ldots,M\}$. 
To judge the experimental characteristics of the proposed tests, we obtained Monte Carlo estimators, say  $\mathrm{PowA}$ and $\mathrm{Pow}$ of the following quantities:
$\mathrm{Pr}_1 \left\{T\left(D_{km}\right)>C_{\alpha}\right\}$ and  $\mathrm{Pr}_1 \left\{T\left(D_{km}\right)>q\left(D_{0m}\right)\right\}$, where  $C_{\alpha}$ is the $\alpha$-level critical value related to the asymptotic $H_0$-distribution of $T$ and  $q\left(D_{0m}\right)$ means a value of the $100(1-\alpha)\%$-quantile of  $T\left(D_{0m}\right)$'s distribution, respectively. 
The criterion $\mathrm{PowA}$ calculated under $D_{0m}$,  $m\in\{1,\ldots,M\}$, examines our current ability to control the TIE rate of a $T$-based test using an approximate $H_0$-distribution of $T$. In this framework,  $\mathrm{PowA}$ calculated under $D_{1m}$,  $m\in\{1,\ldots,M\}$,
displays the expected  power of $T$.  Values of $\mathrm{Pow}$ can be used to evaluate   the actual power levels of $T$, supposing we can accurately control the TIE rates of the corresponding $T$-based test. It can be theoretically assumed that we can correct $T$ to produce a statistic, say $T'$, in order to minimize the distance 
$\left|\mathrm{Pr}_0 \left\{T'\left(D_{0m}\right)>C_{\alpha}\right\}-\alpha\right|$, by employing a method based on, for example, a Bartlett type correction, location adjustments, and/or bootstrap techniques. In this framework, an accurate higher order approximation to the $H_0$-distribution of $T\left(D_{0m}\right)$ might be needed.
In several situations, $\mathrm{Pow}$ could indicate potential abilities to improve practical implementations of  studied tests.  
\subsection{One Sample t-Test for the Mean}\label{sc41} 
In order to examine the $T_N$-based test generated by modifying the t-test, $T_o$,  in Section~\ref{sc31}, the following designs of underlying data distributions were applied:
$D_{k1}:$ $X_1\sim N\left(0.1k,1\right)$;
$D_{k2}:$ $X_i=1-\eta_i+0.1k$ with $\eta_i\sim Exp(1)$;
$D_{k3}:$ $X_i=\eta_i-1+0.1k$;
$D_{k4}:$ $X_i=\left(\xi_i-2\right)/2+0.2k$ with $\xi_i\sim \mathrm{Weibull}(1,2)$;
where $k\in\{0,1\}$ and $i\in\{1,\ldots,n\}$. 
The experimental results presented in Table~\ref{tab1} are the
power comparisons of the  t-test based on $T_o$,  its modification based on $T_N$ and the bootstrap test  $T_B$, the bootstrapped
 version of the t-test, when the significance level, \( \alpha \), of the tests  was supposed to be fixed at 5$\%$.
\begin{table}
	\caption{ Monte Carlo rate of rejections  at $\alpha=0.05$ of the following  statistics: the t-test statistic $T_o$  and its modification $T_N$, defined in Section \ref{sc31}; the t-test statistic's bootstrapped version $T_B$.  \label{tab1}}
	\centering
	\vskip-0.3cm
	\rule[0pt]{17cm}{0.4pt}
	\begin{tabular}{cccccccccccc}
		&   &	\multicolumn{2}{c}{$n=150$}& 	\multicolumn{2}{c}{$n=200$}&	\multicolumn{2}{c}{$n=250$}&	\multicolumn{2}{c}{$n=300$}&	\multicolumn{2}{c}{$n=350$}\\	
		Design&	Test&	PowA & Pow & PowA & Pow& PowA &Pow & PowA &Pow & PowA &Pow\\	
		\vspace{6\baselineskip}					
		$D_{01}$&	$T_o$&0.049&0.050&  0.053&0.050&0.049&0.050&0.050&0.050&0.052&0.050	\\
		        &  $T_N$&0.050 &0.050&  0.054&0.050&0.050&0.050&0.051&0.050&0.053&0.050	\\
		        &  $T_B$&0.054 &0.054&  0.052&0.052&0.050&0.050&0.051&0.051&0.052&0.052\\
		$D_{11}$&	$T_o$&0.337&0.334&	0.409&0.401&0.475&0.478&0.537&0.537&0.588&0.581	\\
		        &  $T_N$&0.337&0.333&	0.412&0.399&0.477&0.476&0.539&0.533&0.590&0.579	\\
		        & $T_B$&0.332&0.332&0.416&0.416&0.463&0.463&0.539&0.539&0.586&0.586\\		     
		$D_{02}$&	$T_o$&0.044 &0.050& 0.046&0.050&0.046&0.050&0.046&0.050&0.044 &0.050\\
		        &  $T_N$&0.044 &0.050& 0.043&0.050&0.044&0.050&0.043&0.050&0.043 &0.050	\\
		        & $T_B$ &0.083 &0.083&0.074	&0.074&	0.068&0.068&0.068&0.068&0.069&0.069\\
		$D_{12}$&	$T_o$&0.345 &0.363& 0.417&0.429&0.487&0.499&0.541&0.554&0.597&0.614\\
		        &  $T_N$&0.479 &0.503& 0.591&0.615&0.679&0.704&0.748&0.774&0.805&0.826\\
		        & $T_B$ &0.386&0.386 &0.447	&0.447&0.507&0.507&0.559&0.559&0.608&0.608\\
		$D_{03}$&	$T_o$&0.055&0.050&0.055&0.050&0.054&0.050&0.053&0.050&0.052&0.050\\
		        &  $T_N$&0.068&0.050&0.068&0.050&0.067&0.050&0.067&0.050&0.065&0.050\\
		        & $T_B$ &0.032&0.032&0.036&0.036&0.037&0.037&0.039&0.039&0.038&0.038\\
		$D_{13}$&	$T_o$&0.331&0.316&0.398&0.379&0.464&0.449&0.528&0.517&0.584&0.579\\
		        &  $T_N$&0.644&0.580&0.732&0.663&0.801&0.749&0.854&0.810&0.890&0.856\\
		        & $T_B$ &0.307&0.307&0.383&0.383&0.454&0.454&0.521&0.521&0.586&0.586\\
		$D_{04}$&	$T_o$& 0.044&0.050&0.046&0.050&0.045&0.050&0.044&0.050&0.046&0.050\\
		        &  $T_N$& 0.042&0.050&0.043&0.050&0.044&0.050&0.044&0.050&0.044&0.050\\
		        & $T_B$ &0.077&0.077&0.079&0.079&0.071&0.071&0.069&0.069&0.068&0.068\\
		$D_{14}$&	$T_o$& 0.790&0.804&0.879&0.885&0.929&0.934&0.961&0.965&0.979&0.980\\
		        &  $T_N$ & 0.940&0.950&0.981&0.984&0.994&0.995&0.998&0.999&0.999&1.000\\
		        & $T_B$  &0.781&0.781&0.863&0.863&0.920&0.920&0.954&0.954&0.973&0.973
		        
	\end{tabular}
	\smallskip
	\centering\small
	\rule[0pt]{17cm}{0.4pt}
\end{table}

Designs $D_{k1}$, $k\in\{0,1\}$, exemplify scenarios, where  $T_o$ is MP. In these cases, values of $\mathrm{Pow}$ testify that $T_o$ is   slightly  superior to $T_N$. 
Designs $D_{k2}$, $k\in\{0,1\}$, correspond to  negatively skewed distributions. In these scenarios, $T_N$ is clearly somewhat better than $T_o$, having  approximately  27$\%$$-$30$\%$ power gains as compared with $T_o$. 
Designs $D_{k3}$, $k\in\{0,1\}$, represent  positively skewed distributions. The proposed test $T_N\left(D_{13}\right)$ is about two times more powerful than $T_o\left(D_{13}\right)$. However, we should note that the asymptotic  TIE rate control related to $T_N\left(D_{03}\right)$ suffers from the skewness of the $H_0$-distribution. According to the values of $\mathrm{Pow}$ computed under $D_{13}$, the procedure $T_N\left(D_{k3}\right)$, $k\in\{0,1\}$, will clearly dominate the strategy  $T_o\left(D_{k3}\right)$, $k\in\{0,1\}$, if the TIE rate control related to $T_N$ could be improved. To this end, for example, a~\cite{Chen}-type approach can be suggested to be applied. The present paper does not aim to achieve improvements of test-algorithms for controlling the TIE rate of $T_N$. The computed values of the criterion $\mathrm{Pow}$ shown in Table~\ref{tab1} confirm that the $T_N$-based strategy is reasonable. The results related to $D_{04}$ and $D_{14}$ support the conclusions  above. Although, under $D_{04}$, the corresponding $\mathrm{PowA}$'s values indicate that  
the Monte Carlo asymptotic TIE rates of $T_N$ are smaller than those related to $T_o$, the proposed test is superior to  $T_o$ in both the $\mathrm{PowA}$ and $\mathrm{Pow}$ contexts under $D_{14}$.  

	In Table~\ref{tab1}, we also report the experimental  results related to the Monte Carlo implementations of  the test based on a bootstrapped version of  the $T_o$ statistic, denoted $T_B$, where $X_1,\ldots,X_n$ are resampled with replacement. In these cases,  
	asymptotic approximations for the corresponding TIE rates were not applied. Thus, we denote   
	  the criterion PowA=Pow.  The applied bootstrap strategy required a substantial computational cost. However, we cannot confirm that the $T_B$-based test is significantly superior to the t-test based on $T_o$, under $D_{01},D_{11},D_{02},\ldots,D_{14}$. Moreover, under the designs $D_{02}$ and $D_{04}$, the 
	the bootstrap t-test cannot be suggested to be used. 
    
\subsection{One Sample Test for the Median}\label{sc42} 
To gain some insight into operating characteristics of the test statistic $T_N$ defined in Section~\ref{sc32}, we considered various designs of underlying data distributions corresponding to the hypotheses $H_0:$ $\nu=0$ and $H_1:$ $\nu>0$, where $\nu$ denotes the median of $X_1$'s distribution. To exemplify the results of the conducted Monte Carlo study, we employ the following schemes:
$D_{01}:$  $X_i=\eta_{i}-\xi_{i}$;
$D_{11}:$  $X_i=1-\eta_{i}$;
$D_{02}:$  $X_i\sim N(0,4)$;
$D_{12}:$  $X_i=\exp(0.5)-\zeta_i$; 
where $\eta_{i}\sim Exp(1)$,  $\xi_{i}\sim Exp(1)$, $\zeta_i\sim LN(0,1)$, and $i\in\{1,\ldots,n\}$.
In this study, attending to the statements presented in Section~\ref{sc32},  
the one-sample, one-sided Wilcoxon-Mann-Whitney test, say W, the $T_o$-based test and its modification, the $T_N$-based test, were implemented.  Note that, for the W test,  the criterion PowA=Pow, since
$H_0$-distributions of the Wilcoxon-Mann-Whitney test statistic do not depend on underlying data distributions.  
Table~\ref{tab2} represents the typical  results observed during  the extensive  power evaluations of  W, $T_o$, and $T_N$, when the significance level, \( \alpha \), of the considered tests  was supposed to be fixed at 5$\%$. 
\begin{table}
	\caption{ Monte Carlo rate of rejections  at $\alpha=0.05$ of the following  statistics: the one-sample Wilcoxon-Mann-Whitney test statistic 
		(W), $T_o$ and its modification, $T_N$, defined in Section \ref{sc32}.  \label{tab2}}
	\centering
	\vskip-0.3cm
	\rule[0pt]{12cm}{0.4pt}
	\begin{tabular}{cccccccc}
		&   &	\multicolumn{2}{c}{$n=25$}& \multicolumn{2}{c}{$n=50$}&	\multicolumn{2}{c}{$n=75$}\\	
		Design&	Test&	PowA & Pow & PowA & Pow& PowA &Pow\\	
		\vspace{6\baselineskip}					
		$D_{01}$&	    W&0.048&0.048&0.049&0.049&0.049&0.049\\
		&	$T_o$&0.029&0.050&0.025&0.050&0.031&0.050\\
		&   $T_N$&0.030&0.050&0.027&0.050&0.028&0.050\\
		$D_{11}$&	    W&0.213&0.213&0.319&0.319&0.417&0.417\\
		&	$T_o$&0.484&0.549&0.674&0.751&0.806&0.862\\
		&   $T_N$&0.660&0.725&0.899&0.936&0.973&0.985\\		
		$D_{02}$&	    W&0.048&0.048&0.050&0.050&0.050&0.050\\
		&	$T_o$&0.049&0.050&0.044&0.050&0.045&0.050\\
		&   $T_N$&0.048&0.050&0.042&0.050&0.045&0.050\\
		$D_{12}$&	    W&0.394&0.394&0.609&0.609&0.755&0.755\\
		&	$T_o$&0.734&0.735&0.919&0.927&0.977&0.980\\
		&   $T_N$&0.843&0.849&0.982&0.986&0.998&0.999\\		
	\end{tabular}
	\smallskip
	\centering\small
	\rule[0pt]{12cm}{0.4pt}
\end{table}
For example, in  scenario $\left\{D_{11},\, n=50\right\}$, $T_N$  improves $T_0$ providing  about a 25\% power gain.

Regarding the two-sided $T_N^2$-based test derived in Section~\ref{sc32}, the following outcomes exemplify the corresponding Monte Carlo power evaluations:
PowA =$0.223$, $0.585$, and $0.838$ provided by the two-sided W-test, $T_o^2$-based test and  
$T_N^2$-based test, respectively, when $n=50$ and generated data satisfy $D_{11}$. Note that, in the scenario above, we can employ the method proposed in~\cite{Thomas}. According to \cite{Thomas}, since $T_o^2>0$, and $T_N^2>0$ are $O_p\left(n^k\right)$, where $k=0$  and $k=1$, under $H_0$ and $H_1$, respectively, the test statistics
$T_{o1}^2=-n\log\left(1-T_{o}^2/n\right)$ and $T_{N1}^2=-n\log\left(1-T_{N}^2/n\right)$
 are reasonable to be examined. These monotonic transformations demonstrated  slight PowA increases of 
 approximately $1.4\%$ and  $1.1\%$ for 
 the $T_{o}^2$- and $T_{N}^2$-based strategies, respectively.

\subsection{Test  for the Center of Symmetry}\label{sc43} 
In this section,  we examine  implementations of the proposed $T_N$-modification of the $T_o$-based test developed in Section~\ref{sc33}.  The $T_o$- and $T_1$-based tests as well as the one-sample, one-sided Wilcoxon-Mann-Whitney test (W) were compared with the $T_N$-based test with respect to the setting depicted in Section~\ref{sc33}. To exemplify the results of the conducted numerical study,  the following designs of data $D=\left\{X_1,\ldots,X_n\right\}$ generations were employed: for $k\in\{0,1\}$ and $i\in\{1,\ldots,n\}$, 
$D_{k1}:$ $X_i\sim N(0.1k,1)$, when  $T_o$ can be expected to be superior to $T_1$, $T_N$, and W; 
$D_{k2}:$ $X_i=\eta_i-\xi_i+0.1k$,  where $\eta_i$ and $\xi_i$ are independent $Exp(1)$-distributed random variables, and then $T_1$ can be expected to be superior to $T_o$, $T_N$, and W;
$D_{k3}:$ $X_i=\zeta_i+0.1k$, where $\zeta_i\sim Unif(-1,1)$;
$D_{k4}:$ $X_i=\epsilon_i-0.5+0.1k$, where $\epsilon_i\sim Beta(0.5,0.5)$. 

Table~\ref{tab3} summarizes the computed Monte Carlo outputs across scenarios $D_{kj}$, $k\in\{0,1\}$,  $j\in\{1,\ldots,4\}$, when $n=50, 150$ and  the significance level, \( \alpha \), of the tests  is supposed to be fixed at 5$\%$. It is observed that: under $D_{01}$ and $D_{11}$, $T_o$ and $T_N$ have very similar  behavior;
under $D_{02}$ with $n=50$, $T_N$ does improve $T_o$ in terms of the TIE rate control;  
under  $D_{12}$, the values of the measurement Pow related to $T_N$ and $T_1$ are  close to each other and greater than those of $T_o$; under $D_{kj}$,   $k\in\{0,1\}$,  $j\in\{3, 4\}$, $T_N$ shows the Monte Carlo power characteristics that outperform those of $T_o$, and W. For example, under $D_{14}$ with $n=50$, $T_N$ has approximately $22\%$, $23\%$, and $65\%$ power gains as compared with W, $T_o$, and $T_1$, respectively.

\begin{table}
	\caption{Monte Carlo power levels at $\alpha=0.05$ of the one-sample Wilcoxon-Mann-Whitney test
		(W) as well as the $T_o$, $T_1$, and  $T_N$-based tests defined in Section \ref{sc33}.  \label{tab3}}
	\centering
	\vskip-0.3cm
	\rule[0pt]{17cm}{0.4pt}
	\begin{tabular}{cccccccccccc}
		&   &	\multicolumn{2}{c}{$n=50$}&	\multicolumn{2}{c}{$n=150$} &   & &	\multicolumn{2}{c}{$n=50$}&	\multicolumn{2}{c}{$n=150$}\\	
		Design&	Test&	PowA & Pow & PowA & Pow& Design&	Test&	PowA & Pow & PowA &Pow\\	
		\vspace{6\baselineskip}					
		$D_{01}$&	    W&0.049&0.049&0.052&0.052&
		$D_{03}$&   
		 W&0.050&0.050&0.049&0.049\\
		        &	$T_o$&0.052&0.050&0.053&0.050&        &
		$T_o$&0.054&0.050&0.050&0.050\\
		        &	$T_1$&0.044&0.050&0.047&0.050&        &
		$T_1$&0.062&0.050&0.057&0.050 \\
		        &   $T_N$&0.054&0.050&0.054&0.050&        &
		$T_N$&0.052&0.050&0.048&0.050 \\  
		$D_{11}$&	    W&0.169&0.169&0.321&0.321&
		$D_{13}$&    
		W&0.301&0.301&0.645&0.645\\
		        &	$T_o$&0.183&0.178&0.338&0.329&        &
	    $T_o$&0.328&0.317&0.678&0.676 \\
		        &   $T_1$&0.128&0.140&0.231&0.238&        &
		$T_1$&0.191&0.160&0.345&0.319 \\
		        &   $T_N$&0.189&0.173&0.341&0.328&        &
		$T_N$&0.390&0.381&0.774&0.782\\
        $D_{02}$&	    W&0.051&0.051&0.049&0.049&$D_{04}$&    W&0.050&0.050&0.050&0.050\\
		        &	$T_o$&0.058&0.050&0.049&0.050&        &$T_o$&0.054&0.050&0.053&0.050\\
		        &	$T_1$&0.026&0.050&0.024&0.050&        &$T_1$&0.089&0.050&0.067&0.050\\
		        &   $T_N$&0.044&0.050&0.035&0.050&        &$T_N$&0.046&0.050&0.048&0.050\\
        $D_{12}$&	    W&0.144&0.144&0.280&0.280&$D_{14}$&    W&0.643&0.643&0.965&0.965\\
		        &	$T_o$&0.132&0.119&0.225&0.225&        &$T_o$&0.637&0.617&0.965&0.962\\
		        &   $T_1$&0.089&0.149&0.206&0.315&        &$T_1$&0.292&0.179&0.493&0.431\\				
		        &   $T_N$&0.134&0.148&0.243&0.304&        &$T_N$&0.826&0.836&0.998&0.998
	\end{tabular}
	\smallskip
	\centering\small
	\rule[0pt]{17cm}{0.4pt}
\end{table}

%
%
%
%
%

Based on the conducted Monte Carlo study, we conclude that the proposed testing strategies exhibit high and stable power characteristics under various designs of alternatives. 

\section{Real Data Example}\label{data}
By blocking the blood flow of the heart, blood clots commonly cause myocardial infarction (MI) events that   lead to heart muscle injury. Heart disease is a leading cause of death affecting about or higher than 20\% of populations regardless of different ethnicities according to the Centers for Disease Control and Prevention, e.g., \cite{DataMI}.

The application of the proposed approach is illustrated by employing a sample from a study that evaluates biomarkers associated with MI. The study was focused on the residents of Erie and Niagara counties, 35$-$79 years of age. The New York State department of Motor Vehicles drivers' license rolls was used as the sampling frame for adults between the age of 35 and 65 years, while the elderly sample (age 65$-$79) was randomly chosen from the Health Care Financing Administration database. The biomarkers called "thiobarbituric acid-reactive substances" (TBARS) and "high-density lipoprotein" (HDL) cholesterol are frequently used as discriminant factors between individuals with (MI=1) and without (MI=0) myocardial infarction disease, e.g., \cite{DataMI}. 

\begin{figure}
	\begin{center}
	\includegraphics[width=3in]{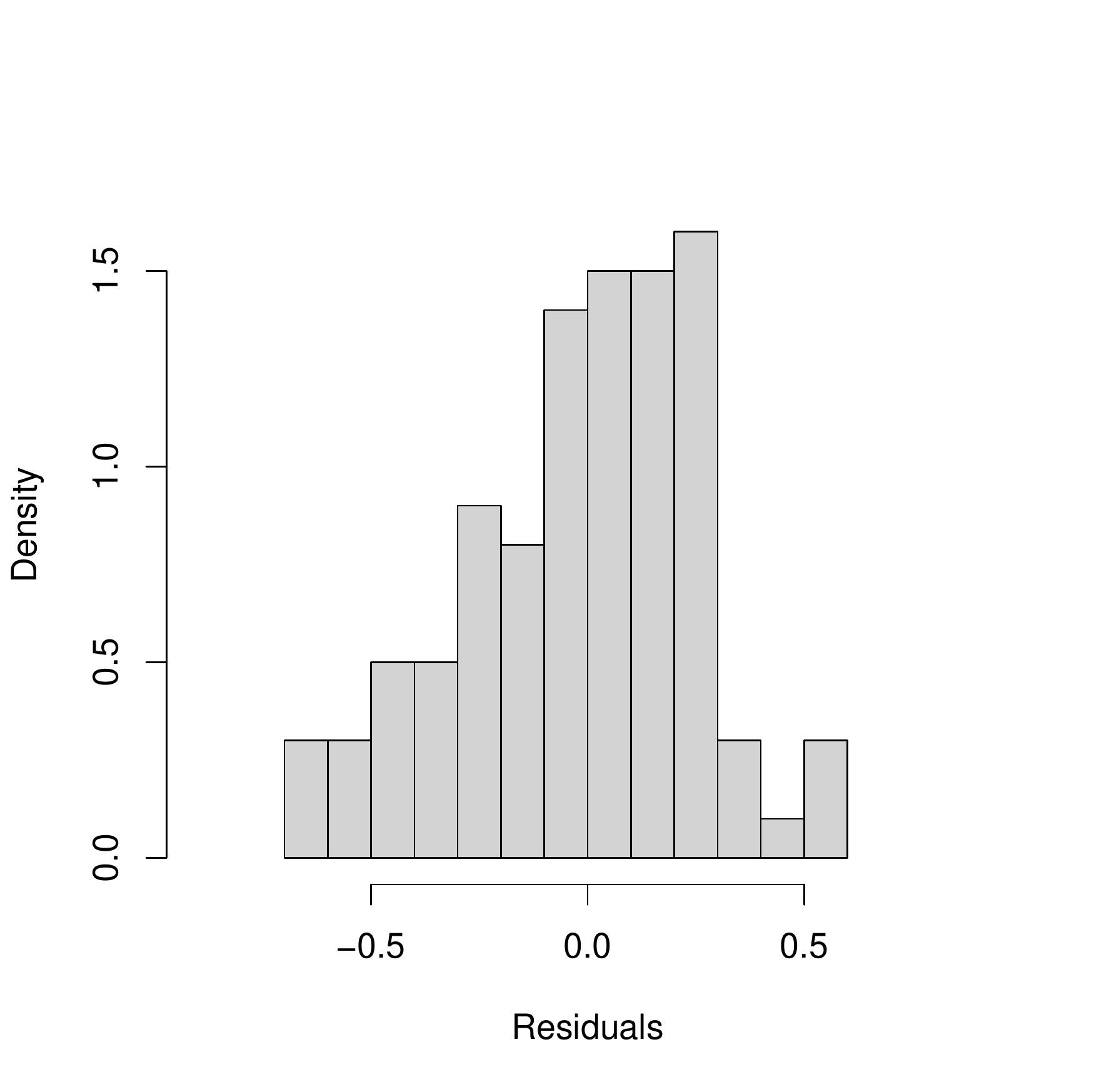}
	\caption{Data-based histogram related to regression residuals $\epsilon_1,\ldots,\epsilon_{n}$.}
	\label{hist}
		\end{center}
\end{figure}

The sample of $2,910$ biomarkers' values was used to estimate the parameters $a$ and $b$ in the linear regression model $Y_i=a+bZ_i+\epsilon_i$ related to \{MI$=1$\}'s cases, where $Y_1,\ldots,Y_{2910}$ are log-transformed HDL-cholesterol measurements, $Z_1,\ldots,Z_{2910}$ denote log-transformed TBARS measurements, and  
$\epsilon_i,\, i\ge 1,$ represent regression residuals with $\text{E}\epsilon_i=0$. It was concluded that $\epsilon_i\simeq Y_i- 4.034-0.045Z_i,\, i\ge 1$ (see Table S1 in the online supplemental material, for details). Assume we aim to investigate the  distribution of $\epsilon_i$ based on $n=100$ biomarkers' values, when $\text{MI}=1$.
In this case, it was observed that the sample mean and variance were $\bar{\epsilon}=\sum_{i=1}^{n}\epsilon_i/n \simeq -0.002$ and 
$\sum_{i=1}^{n}\left(\epsilon_i-\bar{\epsilon}\right)^2/(n-1) \simeq 0.073$, respectively.  
Figure~\ref{hist} depicts the histogram based on corresponding values of $\epsilon_1,\ldots,\epsilon_{100}$. 

In order to test for $H_0:$ $\nu=0$  vs. $H_1:$ $\nu\neq 0$, where $\nu$ is the median of $\epsilon$'s distribution, we implemented the two-sided $T_o^2$-based test and its modification, the $T_N^2$-based test denoted in Section~\ref{sc32}, as well as the  two-sided Wilcoxon-Mann-Whitney test (W). Although  the histogram shown in Figure~\ref{hist} displays a relatively asymmetric  distribution about zero, the $T_o^2$-based test and the W test have demonstrated a  p-value$=0.071$ and p-value$=0.326$, respectively. The proposed $T_N^2$-based test has provided p-value$=0.047$. Then, we organized a Bootstrap/Jackknife type study to examine the power performances of the test statistics. The conducted strategy was that   
a sample with size $n_b<100$  was randomly selected  with replacement from the data $\left\{\epsilon_1,\ldots,\epsilon_{n}\right\}$ to be tested for $H_0$ at a 5\% level of significance.  This strategy was repeated $10,000$ times to calculate the frequencies of the events \{$T_o^2$ rejects $H_0$\}, \{W rejects $H_0$\}, and \{$T_N^2$ rejects $H_0$\}. The obtained experimental powers of $T_o^2$, W, and $T_N^2$ were: $0.238$, $0.106$, $0.535$, when $n_b=90$;
$0.198$, $0.104$, $0.463$, when $n_b=80$;   $0.176$, $0.100$,  $0.415$, when $n_b=70$, respectively. 
The experimental power levels of the tests increase as the sample size $n_b$ increases.  
This study experimentally indicates that the $T_N^2$-based test  outperforms the classical procedures in terms of the power properties when evaluating whether the residuals of the association $Y_i=a+bZ_i+\epsilon_i,\, i\ge1$, are distributed asymmetrically about zero. That is, the proposed test can be expected to be more sensitive as compared with the known methods to rejecting the null hypothesis $H_0:$ $\nu=0$ vs. $H_1:$ $\nu\neq 0$, in this study.   

\section{Concluding Remarks}\label{sc5}
The present article has provided a theoretical framework for evaluating and constructing powerful data-based tests.  
The contributions in this article have touched on the principles of characterizing most powerful statistical decision-making mechanisms. Proposition~\ref{pr3} provides a method
for one-to-one mapping the term "most powerful" to the properties of  test statistics' distribution functions via analyzing the behavior of corresponding likelihood ratios. 
We  demonstrated that the derived characterization of MP tests can be associated with  a  principle of sufficiency.
The concepts shown in Section~\ref{sc2} have been  applied to  improving test procedures 
by accounting for the relevant ancillary statistics.
 Applications of the presented  theoretical framework have been employed to display efficient modifications of the one-sample t-test,  the test for the median, and the test for the center of symmetry, in  nonparametric settings. 
The effectiveness of the proposed nonparametric decision-making procedures  in maintaining relatively high power has been confirmed using simulations and a real data example across various scenarios based on  samples from relatively skewed distributions. 
We also note the following remarks. (a) Propositions~\ref{pr4} and \ref{pr6}  can be applied to  different decision-making problems. 
(b)  
Effective corrections of the classical t-test can be of  theoretical and applied interest. 
The modification of the t-test  as per Section~\ref{sc31} involves using the estimation of the  third central moment. This moment plays a role in some corrections of the  t-test structure for adjusting its  null distribution when the underlying data are asymmetric, e.g.,~\cite{Chen}. Overall, the proposed modification is somewhat different from those that are used to improve control of the TIE rates of  t-test type procedures. Thus, in general, basic ingredients of the methods mentioned above can be combined.  
(c) The scheme presented in Section~\ref{sc32} can be easily revised to develop a test for quantiles, by using the observation that $\bar{X}$ and the sample $p$th quantile, say $X_{(pn)}$, are asymptotically bivariate normal with  
\begin{eqnarray*}
	&&\mathrm{cov}\left(\bar{X},X_{(pn)}\right)\,=\,f\left(\nu_p\right)^{-1}
	\mathrm{E}\left(X_1-\nu_p\right)\left\{pI\left(X_1>\nu_p\right)-(1-p)I\left(X_1\le\nu_p\right)\right\},
	\\
	&&\mathrm{var}\left(X_{(pn)}\right)\,=\,f\left(\nu_p\right)^{-2}p(1-p),\,\,\, \mathrm{Pr}\left(X_1<\nu_p\right)\,=\,p.
\end{eqnarray*} 
(d) Sections~\ref{sc3}$-$\ref{data}  have exemplified applications of the treated MP principle in the nonparametric settings.
In many parametric problems, the corresponding likelihood ratios do not  have explicit forms or have very complicated shapes, e.g., when testing statements are based on longitudinal data, dependent observations, multivariate outcomes, data subject to different sorts of errors, and/or missing-values mechanisms.
In such cases, issues related to comparing/developing tests via the considered MP principle  can be employed.

A plethora of decision-making algorithms touches on most fields of statistical practice. Thus, it is not practical in one paper to focus on all the relevant theory and examples. This paper studies only one approach to characterize  a class  of MP mechanisms. That is, there are many potential future directions that seem to be promising targets for research, including, for example: (i) examinations of the relationships
between general MP characterizations, Basu's theorem-type results (e.g.,~\citealp{Basu}),  the concepts of sufficiency, completeness, and ancillarity   under   different statements of decision-making policies. In this aspect, for example,  a research question can be as follows: when can we claim that a statistic $T$ is MP iff $T$ and $A$ are independently distributed, for any ancillary statistic $A$? 
(ii) Various parametric and nonparametric applications of MP characterizations in different settings can be developed.
(iii) Proposition~\ref{pr4} can be used and extended  to compare different statistical procedures in practice.
(iv) In  light of the present MP principle, relevant evaluations of optimal combinations of test statistics (e.g.,~\citealp{Combi}) can be proposed. 
(v)  Large sample properties of  test statistics  modified with respect to  MP characterization can be analyzed.
(vi) Perhaps,  Proposition~\ref{pr6} can be integrated into various testing developments, where  characterizations of underlying data distributions under corresponding hypotheses can be used to define relevant ancillary statistics.  
Leaving these topics to the future, it is hoped that the present paper will convince the readers of  the benefits of studying different aspects and characterizations related to MP data-based decision-making techniques.      

\section*{Supplementary Materials}
The online supplementary materials contain: the proofs of the theoretical results presented in the article; and Table S1 that displays the analysis of variance related to the linear regression fitted to the observed log-transformed HDL-cholesterol measurements using the log-transformed TBARS measurements as a factor, in Section~\ref{data}. 

\bibliographystyle{Chicago}

\bibliography{paperMPtest.bib}

\begin{center}	
\Large{ \bf Supplemental Materials to "{\em A Characterization of Most(More) Powerful Test Statistics with Simple Nonparametric Applications}"}
\end{center}	
\spacingset{1.45} 

\section*{APPENDIX: PROOFS}
This Appendix comprises the necessary proofs  to establish the propositions presented in the present paper.
\subsection*{Proof of Proposition 1.1}
Consider 
\begin{eqnarray*}
	&&{\Pr}_{1} \left(u-s\le \Lambda\le u\right)=\text{E}_{1 }
	I\left(u-s\le \Lambda\le u\right)
	\,=\, \int I\left(u-s\le \Lambda\le u\right)f_{1} 
	\\
	&&\quad =\,\int I \left(u-s\le \Lambda\le u\right)\frac{f_{1 } }{f_{0} } f_{0} \,\,\,=\,\int I \left(u-s\le \Lambda\le u\right)\, \Lambda\,\,f_{0},
\end{eqnarray*}
where $I(.)$ means the  indicator function and $u,s$ are not random variables. This implies the inequalities
\begin{eqnarray*}
	{\Pr}_{1} \left(u-s\le \Lambda\le u\right)&\le&
	\int I \left(u-s\le \Lambda\le u\right)u f_{0} \,\,\,\,\,=\,u{\Pr}_{0} \left(u-s\le \Lambda\le u\right),
	\\
	{\Pr}_{1} \left(u-s\le \Lambda\le u\right)&\ge&
	\int I \left(u-s\le \Lambda\le u\right)(u-s) f_{0} \\
	&&=\,(u-s){\Pr}_{0} \left(u-s\le \Lambda\le u\right).	 
\end{eqnarray*}
Dividing these inequalities by $s$ and employing $s\to 0$, we complete the proof.
\subsection*{Proof of Proposition 2.1}
We can write
\begin{eqnarray*}
	f_{1 }^{\Lambda,A} (u,v)&=&\frac{d}{dv}\,\, {\mathop{\lim }\limits_{s\to 0}} \, \frac{1}{s} {\Pr}_{1 } \left(u-s<\Lambda<u,A<v\right),
\end{eqnarray*} 
where 
\begin{eqnarray*}
	&&{\Pr}_{1 } \left\{u-s<\Lambda(D)<u,A(D)<v\right\}=\text{E}_{1 } I\left\{u-s<\Lambda(D)<u,A(D)<v\right\}\\ 
	&&\quad =\int I\left\{u-s<\Lambda(D)<u,A(D)<v\right\}\, f_{1 } (D)
	\\
	&&\quad=\int I\left\{u-s<\Lambda(D)<u,A(D)<v\right\}\, \Lambda(D)\,\, f_{0} (D).
\end{eqnarray*} 
Noting that 
\begin{eqnarray*}
	\int I\left(u-s<\Lambda<u,A<v\right)\, \Lambda\,\, f_{0}(D)&\le& u\int I\left(u-s<\Lambda<u,A<v\right)\, f_{0} (D), \\
	\int I\left(u-s<\Lambda<u,A<v\right)\, \Lambda\,\, f_{0} (D)&\ge& \left(u-s\right)\int I\left(u-s<\Lambda<u,A<v\right) \, f_{0} (D),
\end{eqnarray*} 
we have 
\begin{eqnarray*}
	{\Pr}_{1 } \left(u-s<\Lambda<u,A<v\right)&\le& u\, {\Pr}_{0} \left(u-s<\Lambda<u,A<v\right), 
	\\
	{\Pr}_{1 } \left(u-s<\Lambda<u,A<v\right)&\ge& \left(u-s\right)\,{\Pr}_{0} \left(u-s<\Lambda<u,A<v\right).
\end{eqnarray*} 
This completes the proof in a similar manner to that of Proposition 1.1.
\subsection*{Proof of Proposition 2.2}
The proof of Proposition 2.2 consists of the two steps below:  
(1) When $T=\Lambda$, i.e. $T$ is the MP test statistic,  Proposition 2.1 says $f_{1 }^{T,A} (u,v)=u{\rm \; }f_{0}^{T,A} (u,v)$. 
(2) When $f_{1 }^{T,A} (u,v)=u{\rm \; }f_{0}^{T,A} (u,v)$ is satisfied, we consider the elementary inequality: for all $B$ and $C$,
\[(B-C)\left\{I(B\ge C)-\delta \right\}\ge 0,\] 
where $I\left(.\right)$ is the indicator function and $\delta \in \left[0,1\right]$. Let $\delta =I\left(\Lambda\ge U\right)$ denote the $H_0$-rejection rule of the likelihood ratio test with the fixed threshold $U$, such that: we reject $H_0$, if $\delta=1$. Letting $B=T$, we have 
\[\text{E}_{0} \left\{\left(T-C\right)I\left(T\ge C\right)\right\}\ge \text{E}_{0} \left\{\left(T-C\right)I\left(\Lambda\ge U\right)\right\}.\] 
Thus, fixing the thresholds $C,\, U$ such that the TIE rates ${\Pr}_{0} \left(T\ge C\right)={\Pr}_{0} \left(\Lambda\ge U\right)$, we obtain
$$\text{E}_{0} I\left(T\ge C\right)T\,\ge\, \text{E}_{0} I\left(\Lambda\ge U\right)T,$$ 
where
\[\text{E}_{0} I\left(\Lambda\ge U\right)T\,=\,\text{E}_{0} I\left(\Lambda\ge U\right)f_{1 }^{T,\Lambda} (T,\Lambda)/f_{0}^{T,\Lambda} (T,\Lambda) ,\] 
since $u=f_{1 }^{T,A} (u,v)/f_{0}^{T,A} (u,v)$,  for any statistic $A$ and all $u\ge0$. That is to say,
\begin{eqnarray*}
	&&\text{E}_{0} I\left(\Lambda\ge U\right)T={\int\!\!\!\!\int}I\left(v\ge U\right)u\,f_{0}^{T,\Lambda} (u,v)dudv 
	\\
	&&\quad={\int\!\!\!\!\int}I\left(v\ge U\right)\frac{f_{1 }^{T,\Lambda} (u,v)}{f_{0}^{T,\Lambda} (u,v)} f_{0}^{T,\Lambda} (u,v)dudv  
	={\int\!\!\!\!\int}I\left(v\ge U\right)f_{1 }^{T,\Lambda} (u,v)dudv 
	\\
	&&\quad=\int I\left(v\ge U\right)\int f_{1 }^{T,\Lambda} (u,v)du dv 
	=\int I\left(v\ge U\right)f_{1 }^{\Lambda} (v)dv \,\,=\,\, {\Pr}_{1 } \left(\Lambda\ge U\right). 
\end{eqnarray*}
Thus,
\[\text{E}_{0} I\left(T\ge C\right)T\ge {\Pr}_{1 } \left(\Lambda\ge U\right),\] 
where 
\begin{eqnarray*}
	\text{E}_{0}  I\left(T\ge C\right) T &=&\int I\left(u\ge C\right)uf_{0}^{T} (u) du
	=\int I\left(u\ge C\right)
	\frac{f_{1}^{T} (u)}{f_{0}^{T} (u)} f_{0}^{T} (u) du
	\\ 
	&=&{\Pr}_{1} \left(T\ge C\right), 
\end{eqnarray*}
since \textbf{$f_{1 }^{T,A} (u,v)=u{\rm \; }f_{0}^{T,A} (u,v)$ }implies that $\int f_{1 }^{T,A} (u,v) dv=u{\rm \; }\int f_{0}^{T,A} (u,v) dv$, i.e.

\noindent  $f_{1 }^{T} (u)/f_{0}^{T} (u)=u$. This provides 
$${\Pr}_{1 } \left(T\ge C\right)\ge {\Pr}_{1 } \left(\Lambda\ge U\right),$$
meaning that if $f_{1 }^{T,A} (u,v)=u f_{0}^{T,A} (u,v)$ then $T$ is MP, since $\Lambda$ is MP.
The proof is complete. 

\subsection*{Proof of Proposition 2.3}
It is clear that $T=\Lambda=f_1(D)/f_0(D)$, i.e. $T$ is MP, implies 	$\mathrm{E}_{1 }\left\{g(D)\right\}=\mathrm{E}_{0 }\left\{g(D)T(D)\right\}$.

Now, we assume that $T$ satisfies 	$\mathrm{E}_{1 }\left\{g(D)\right\}=\mathrm{E}_{0 }\left\{g(D)T(D)\right\}$ and the event $T\ge C$ suggests to reject $H_0$. Consider the elementary inequality: for all $B$ and $C$,
\[(B-C)\left\{I(B\ge C)-\delta \right\}\ge 0,\,\,\, \delta \in \left[0,1\right].\]
Suppose 
$\delta =I\left(\Lambda \ge U\right)$ denotes the $H_0$-rejection rule of the likelihood ratio test with the fixed threshold $U$, such that: we reject $H_0$, if $\delta=1$ and ${\Pr}_{0} \left(T\ge C\right)={\Pr}_{0} \left(\Lambda\ge U\right)=\alpha$. In the elementary inequality above, we define  $B=T$ and then obtain
\begin{eqnarray*}
	0\le \mathrm{E}_{0 } (T-C)\left\{I(T\ge C)-I\left(\Lambda \ge U\right) \right\}\,=\, \mathrm{E}_{0 } \left\{I(T\ge C) T\right\}-\mathrm{E}_{0 }\left\{ I\left(\Lambda \ge U\right) T\right\}.
\end{eqnarray*}
This yields that 
\begin{eqnarray*}
	0\le  \mathrm{E}_{1 } I(T\ge C) -\mathrm{E}_{1 } I\left(\Lambda \ge U\right),
\end{eqnarray*} 
since $\mathrm{E}_{0 }\left\{g(D)T(D)\right\}=\mathrm{E}_{1 }\left\{g(D)\right\}$, for every $g\in [0,1]$.
But, ${\Pr}_1\left(\Lambda \ge U\right)\ge {\Pr}_1\left(T \ge C\right)$, since $\Lambda$ is MP. Therefore ${\Pr}_1\left(\Lambda \ge U\right)= {\Pr}_1\left(T \ge C\right)$, meaning that $T$ is MP.
The proof is complete. 
\subsection*{Proof of Proposition 2.4}
To prove Proposition 2.4, we note that, as mentioned above (see, e.g., the proof of Proposition 2.2), $f_{1}^{T_{1} ,T_{2} } (u,v) =uf_{0}^{T_{1},T_{2} }(u,v)$ implies  $f_{1}^{T_{1} } (u) =uf_{0}^{T_{1}}(u)$
that leads to  $\text{E}_{0} I\left(T_1\ge C\right)T_1={\Pr}_{1 } \left(T_1\ge C\right)$, for a fixed threshold $C$. Since $f_{1}^{T_{1} ,T_{2} } (u,v) =uf_{0}^{T_{1},T_{2} }(u,v)$, we also have that 
\begin{eqnarray*}
	&&\text{E}_{0} I\left(T_2\ge U\right)T_1={\int\!\!\!\!\int}I\left(v\ge U\right)u\,f_{0}^{T_1,T_2} (u,v)dudv 
	\,\,=\,\, {\Pr}_{1 } \left(T_2\ge U\right), 
\end{eqnarray*}
where a fixed  $U$ satisfies ${\Pr}_{0} \left(T_1\ge C\right)$ $={\Pr}_{0} \left(T_2\ge U\right)=\alpha$, $\alpha$ denotes the TIE rate of the tests. Then, in a similar manner to   the proof  of  Proposition 2.2, we apply the inequality 
$ \text{E}_{0} I\left(T_1\ge C\right)T_1$ 
$\ge \text{E}_{0} I\left(T_2\ge U\right)T_1$ to complete the proof. 
\subsection*{Proof of Proposition 2.5}
By virtue of the statements presented above Proposition 2.5 in the article, the proof is straightforward.
\subsection*{Proof of Proposition 3.1}
Consider the bivariate characteristic function
\begin{eqnarray*}
	&&\text{E}_{1}e^{\mathrm{i}t_1T_N+\mathrm{i}t_2T}\,=\,
	\text{E}_{1}e^{\mathrm{i}t_1T_N+\mathrm{i}t_2\psi\left(T_N,A\right)}\,=\,
	\int\int e^{\mathrm{i}t_1u+\mathrm{i}t_2\psi\left(u,x\right)} f_1^{T_N,A}(u,x)dudx
	\\
	&& \qquad =\, \int\int e^{\mathrm{i}t_1u+\mathrm{i}t_2\psi\left(u,x\right)} f_1^{T_N}(u) f_1^{A}(x)dudx
	\,=\,\int\int e^{\mathrm{i}t_1u+\mathrm{i}t_2\psi\left(u,x\right)} f_1^{T_N}(u) f_0^{A}(x)dudx
	\\
	&& \qquad =\, \int\int L^{T_N}(u)\, e^{\mathrm{i}t_1u+\mathrm{i}t_2\psi\left(u,x\right)} f_0^{T_N}(u) f_0^{A}(x)dudx
	\,=\, \text{E}_0\left\{ L^{T_N}(T_N)\, e^{\mathrm{i}t_1 T_N+\mathrm{i}t_2\psi\left(T_N,A\right)}\right\}
	\\
	&&=\text{E}_0\left\{ L^{T_N}(T_N)\, e^{\mathrm{i}t_1 T_N+\mathrm{i}t_2T}\right\},
\end{eqnarray*}
where \( {{\mathrm{i}}^{2}}=-1\) and \(t_j\in {{\mathbb{R}}^{1}} \), $j\in\{1,2\}$.
In this equation, the  Fourier  transformations of $f_{1}^{T_{N} ,T} (u,v)$ and $L^{T_N}(u)\,f_{0}^{T_{N} ,T} (u,v)$ are equivalent  that means  $f_{1}^{T_{N} ,T} (u,v)$  $=L^{T_N}(u)\,f_{0}^{T_{N} ,T} (u,v)$. As mentioned above Proposition 3.1 in the article, since the  ratio  
$L^{T_N}(u)$ is a monotonically increasing function, we can  without loss of generality assume that 
$L^{T_N}(T_N)=T_N$. Then, we have  
\begin{eqnarray*}
	\text{E}_{1}e^{\mathrm{i}t_1T_N+\mathrm{i}t_2T}&=&\text{E}_0\left\{ T_N \, e^{\mathrm{i}t_1 T_N+\mathrm{i}t_2T}\right\},
\end{eqnarray*}
yielding  $f_{1}^{T_{N} ,T} (u,v)$  $=u\,f_{0}^{T_{N} ,T} (u,v)$. Now, by virtue of Proposition 2.4,  the proof of Proposition 3.1 is complete, when $L^{Y}(u)$ is a monotonically increasing function. Similarly, we can consider the case, where 
$L^{Y}(u)$ is a monotonically decreasing function. The proof is complete. 

\section*{TABLE S1}
\begin{table}
	\caption{The analysis of variance related to the linear regression fitted to the observed log-transformed HDL-cholesterol measurements, $Y_i,i=1,\ldots,2910$, using the log-transformed TBARS measurements, $Z_i,i=1,\ldots,2910$, as a factor, in Section 5.}
	\centering
	\vskip-0.3cm
	\rule[0pt]{17cm}{0.4pt}
	\begin{tabular}{ccccc}
		\multicolumn{2}{l} {Residuals:}& & & \\
		Min  &       1Q &             Median &        3Q &          Max\\ 
		-0.80510&     -0.15987&      0.03547&     0.20073&     0.60918 \\
		\multicolumn{2}{l} {Coefficients:} & & & \\
		&Estimate&                  Std. Error&                t value&     $\Pr(>|t|)$\\    
		(Intercept)&       4.0342&  0.03131&    207.371&          $<$ 2e-16 \\
		$Z$        &      -0.0448&  0.00568&    -7.786&          9.55e-15\\
		\multicolumn{4}{l}{Residual standard error: 0.2693 on 2908 degrees of freedom}& \\
		\multicolumn{4}{l}{Multiple R-squared:  0.02042,   Adjusted R-squared:  0.02008} & \\
		\multicolumn{4}{l}{F-statistic: 60.62 on 1 and 2908 DF, p-value: 9.552e-15} & 		
	\end{tabular}
	\smallskip
	\centering\small
	\rule[0pt]{17cm}{0.4pt}
\end{table}



\end{document}